\documentclass[12pt,a4paper]{article}

\usepackage{graphicx}
\usepackage{color}
\usepackage{amsmath}
\usepackage{amssymb}
\usepackage{amscd}
\usepackage{amsthm}
\usepackage{delarray}
\usepackage{enumerate}
\usepackage[T1]{fontenc}
\usepackage{inputenc}
\usepackage{enumerate}
\input xy

\begin{document}

\xyoption{all}
\setlength{\parindent}{5mm}
\renewcommand{\leq}{\leqslant}
\renewcommand{\geq}{\geqslant}
\newcommand{\N}{\mathbb{N}}
\newcommand{\sph}{\mathbb{S}}
\newcommand{\Z}{\mathbb{Z}}
\newcommand{\R}{\mathbb{R}}
\newcommand{\C}{\mathbb{C}}
\newcommand{\F}{\mathbb{F}}
\newcommand{\K}{\mathbb{K}}
\newcommand{\RN}{\mathbb{R}^{2n}}
\newcommand{\derive}[2]{\frac{\partial{#1}}{\partial{#2}}}
\renewcommand{\S}{\mathbb{S}}
\renewcommand{\H}{\mathbb{H}}
\newcommand{\eps}{\varepsilon}
\theoremstyle{plain}
\newtheorem{theo}{Theorem}[section]
\newtheorem{prop}[theo]{Proposition}
\newtheorem{lem}[theo]{Lemma}
\newtheorem{definit}[theo]{Definition}
\newtheorem{corol}[theo]{Corollary}

\title{On some completions of the space of Hamiltonian maps}
\author{\textsc{Vincent Humili{\`e}re}}
\date{}
\maketitle
\normalsize
\begin{center}
Centre de Mathématiques Laurent Schwartz

UMR 7640 du CNRS

Ecole Polytechnique - 91128 Palaiseau, France

\texttt{vincent.humiliere@math.polytechnique.fr}
\end{center}

\begin{abstract}

In his paper \cite{V1}, C. Viterbo defined a distance on the set of
Hamiltonian diffeomorphisms of $\R^{2n}$ endowed with the standard
symplectic form $\omega_0=dp\wedge dq$. We study the completions of
this space for the topology induced by Viterbo's distance and some
others derived from it, we study their different inclusions and give
some of their properties.

In particular, we give a convergence criterion for these distances
that allows us to prove that the completions contain non-ordinary
elements, as for example, discontinuous Hamiltonians. We also prove
that some dynamical aspects of Hamiltonian systems are preserved in
the completions.

\end{abstract}

\section{Introduction.}

Given an open subset $U$ in $\RN$, we denote by $Ham(U)$ the set of all 1-periodic time dependent Hamiltonian functions
$\R\times\RN\to\R$ whose support for fixed time is compact and contained in $U$. We will write $Ham$ for $Ham(\RN)$.

Given a Hamiltonian function $H\in Ham$, its symplectic gradient
(i.e the unique vector field $X_H$ satisfying
$dH=\iota_{X_H}\omega_0$) generates a Hamiltonian isotopy
$\{\phi_H^t\}$. The set of Hamiltonian diffeomorphisms generated by
an element $H$ in $Ham(U)$ will be denoted by
$\mathcal{H}(U)=\{\phi_H=\phi_H^1\,|\,H\in Ham(U)\}$, and we will
write $\mathcal{H}$ for $\mathcal{H}(\RN)$. Finally, we call
$\mathcal{L}=\{\phi(0_n)\,|\,\phi\in\mathcal{H}\}$, the set of
Lagrangian submanifolds obtained from the zero section $0_n\subset
T^*\R^n=\RN$, by a Hamiltonian isotopy with compact support.

As usual, we denote Viterbo's distance on $\mathcal{L}$ or $\mathcal{H}$ by $\gamma$ (see \cite{V1}). Convergence with respect to $\gamma$ is called c-convergence.

\medskip
Our main goals in this paper is to understand the completion
$\overline{\mathcal{H}}^{\gamma}$ of the metric space
$(\mathcal{H},\gamma)$, to give some convergence criterion (section
\ref{critere}) and to compare it with the convergence for Hofer's
distance $d_H$ (see \cite{HZ}, chapter 5 section 1).

 The notion of  $C^0$ symplectic topology has been studied by many
authors, starting from the work of Eliashberg and Gromov on the
$C^0$ closure of the group of symplectic diffeomorphisms, to the
later results of Viterbo (\cite{V1}) and Hofer (\cite{H1}).

More recently Oh (\cite{OH}) gave a deep study of several versions
of $C^0$ Hamiltonians. However, our definition seems to differ from
his, since in all his definitions, he needs the Hamiltonians to be
continuous, while our study starts as we drop this assumption.

\medskip
Let us now state our main results. For convenience, they will be
restated throughout the paper. In section \ref{critere}, we
introduce a symplectic invariant $\xi_{\infty}$ associated to any
subset of $\RN$, and prove that

\begin{theo}\label{th intro}
Let $(H_k)$ be a sequence of Hamiltonians in $Ham$, whose supports
are contained in a fixed compact set. Suppose there exist a
Hamiltonian $H\in Ham$ and a compact set $K\in\RN$ with
$\xi_{\infty}(K)=0$, such that $(H_k)$ converges uniformly to $H$ on
every compact set of $\R\times(\RN -K)$. Then $(\phi_{H_k})$
converges to $\phi_H$ for $\gamma$.
\end{theo}

Examples of sets $K$ with $\xi_{\infty}(K)=0$ are given by compact submanifolds of dimension $d\leq n-2$.

\medskip
Viterbo's distance $\gamma$ is defined on $\mathcal{H}$, but we can
define for any $H,K\in Ham$
$$\gamma_u(H,K)=\sup\{\gamma(\phi^t_{H},\phi^t_K)\,|\,t\in [0,1]\}, $$
to get a new distance on $Ham$ (we give several variants of this
definition). Then the following proposition allows to extend the
notion of Hamiltonian flow.

\begin{prop} If we consider the respective completions
$\overline{\mathcal{H}}^{\gamma}$ and $\overline{Ham}^{\gamma_u}$ of
the metric spaces $(\mathcal{H},\gamma)$ and $(Ham, \gamma_u)$, then
the map $(H,t)\mapsto\phi_H^t$, $Ham\times\R\to\mathcal{H}$ induces
 a map
$\overline{Ham}^{\gamma_u}\times\R\to\overline{\mathcal{H}}^{\gamma}$.

The induced map associates to any element $H$ in
$\overline{Ham}^{\gamma_u}$ a path in
$\overline{\mathcal{H}}^{\gamma}$ that we will call the
\textnormal{generalized Hamiltonian flow} generated by $H$.
\end{prop}

We then show that some aspects of Hamiltonian dynamics can be
extended to the completions (section \ref{extension dynamics}): We
can define a natural action of a generalized flow on a Lagrangian
submanifold. We can also associate to it a support and extend the
notion of first integral.

To some of them, it is also possible, as we prove in section
\ref{Hamilton Jacobi}, to associate a solution to the
Hamilton-Jacobi equation:
 $$\derive{u}{t}+H\left(t,x,\derive{u}{x}\right)=0.$$
Indeed, a $\gamma_2$-Cauchy sequence of Hamiltonians gives a
$C^0$-Cauchy sequence of solutions (where $\gamma_2$ denotes one
variant of the distance $\gamma_u$ we mentioned above).

In section we give examples of elements in both completions
$\overline{\mathcal{H}}^{\gamma}$ and $\overline{Ham}^{\gamma_u}$
that can be described in a much more concrete way than their
abstract definition (as equivalence classes of Cauchy sequences).
More precisely, we prove
\begin{prop} There is a one-one map
$$\mathfrak{F}^{\infty}\to\overline{Ham}^{\gamma_u},$$
where $\mathfrak{F}^{\infty}$ denotes the set of all functions
$H:\R\times\RN\to\R\cup\{+\infty\}$ such that:
\begin{description}
\item[(i)] $H$ is continuous on $\R\times\RN$,
\item[(ii)] $H$ vanishes at infinity: $\forall\eps>0,\exists r,(|x|>r\Rightarrow (|\forall t, H(t,x)|<\eps))$,
\item[(iii)] there exists a zero capacity set (e.g. an infinitesimally displaceable set), that
contains all the points $x$ where $H(t,x)$ is $+\infty$ for time $t$
\item[(iv)] $H$ is smooth on $\R\times\RN-H^{-1}(\{+\infty\})$.
\end{description}
\end{prop}

\medskip
Finally, let us mention that although we developed our theory on
$\RN$, we can  reasonably expect similar results (except those of
sections \ref{section action lag} and \ref{Hamilton Jacobi}) on any
compact symplectic manifold satisfying
$$\omega|_{\pi_2(M)}=0\text{ and } c_1|_{\pi_2(M)}=0.$$
Indeed, on these manifolds, Schwarz defined in \cite{Sc} a distance
which is entirely analogous to Viterbo's.

\medskip
\noindent \textbf{Organization of the paper. }In Section
\ref{invariants} we give the definitions of the objects used in the
paper. For the reader's convenience, we first recall the
construction of Viterbo's distance $\gamma$ (\ref{subsection
invariants de viterbo}) which is based on the theory of generating
functions for Lagrangian submanifolds (\ref{subsection gfqi}). We
also remind the reader of the different symplectic capacities
constructed from $\gamma$ (\ref{def cap}). Finally we introduce our
new distances derived from $\gamma$ (\ref{distances derivées}).

Section \ref{critere} is fully devoted to the proof of our
convergence criterion. Examples of cases where it holds is then
given in \ref{xi petit}.

In Section \ref{extension dynamics} we define the completions of
$Ham$ and $\mathcal{H}$ and show that some aspects of Hamiltonian
dynamics that can be extended to the completions.

In Section \ref{description elements} we discuss some interesting
examples of elements of the completions.

Our results on the Hamilton-Jacobi equation are given in Section
\ref{Hamilton Jacobi}.

Finally, we prove in Appendix a "reduction inequality" usefull to
prove then all the inequalities between the distances considered in
the paper.

\medskip
\noindent \textbf{Acknowledgments. } I am grateful to my supervisor
C. Viterbo for his advices. I also want to thank my friends M. Affre
and N. Roy for spending hours correcting my awful English.

\section{Symplectic invariants.}\label{invariants}

In this section we give the definitions of all the objects we will
use in the sequel. We first recall the definition of Viterbo's
distance, defined first for Lagrangian submanifolds with the help of
generating functions, and then for Hamiltonian diffeomorphisms (see
\cite{V1}).

\subsection{Generating functions quadratic at infinity.}\label{subsection gfqi}

Let $L$ be a Lagrangian submanifold of the cotangent bundle $T^*M$ of a smooth manifold $M$. We say that $L$ admits a \textit{generating function} if there exists an integer $q>0$ and a smooth function $S:M\times\R^q\to\R$ such that $L$ can be written
$$L=\left\{ (x,p)\in T^*M\,|\,\exists\xi\in\R^q,\derive{S}{\xi}(x,\xi)=0 \text{ and } \derive{S}{x}(x,\xi)=p\right\}.$$
Such function $S$ is called a \textit{generating function quadratic at infinity} (or just ``g.f.q.i'') if there exists a non degenerate quadratic form $Q$ on $\R^q$ and a compact $K\subset M\times\R^q$ such that, $\forall (x,\xi)\notin K, S(x,\xi)=Q(\xi)$.

For instance, any quadratic form on $\R^q$ viewed as a function on
$M\times\R^q$ is a g.f.q.i of the zero section $0_M\subset T^*M$.
J.C. Sikorav proved in \cite{S1} that the property of having a
g.f.q.i is invariant by Hamiltonian isotopy with compact support.
For this reason we will be interested in the set $\mathcal{L}$ of
Lagrangian submanifolds, images of the zero section by a Hamiltonian
isotopy with compact support.

Furthermore, C. Viterbo and D. Théret proved that the g.f.q.i's of a given Lagrangian submanifold are essentially unique. Before stating this result, let us introduce the following definitions: For a given function $S:M\times\R^q\to\R$, we call a \textit{stabilisation} of $S$ any function $S':M\times\R^q\times\R^{q'}\to\R$ of the form $S'(x,\xi,\xi')=S(x,\xi)+q(\xi')$, where $q$ is a non-degenerate quadratic form on $\R^{q'}$. In addition, two functions $S,S':M\times\R^q\to\R$ are said \textit{equivalent} if there exists a diffeomorphism $\phi$ of $M\times\R^q$ and a real $C$ such that $S'=S\circ\phi+C$.

\begin{theo}[\cite{V1,T}]\label{uniqueness}
Suppose $S$, $S'$ are two g.f.q.i's of the same Lagrangian submanifold in $\mathcal{L}$. Then, up to stabilisation, $S$ and $S'$ are equivalent.
\end{theo}

This result allows to associate symplectic invariants to any element of $\mathcal{L}$.

\subsection{Invariants defined by minimax and a distance on the group of Hamiltonian
diffeomorphisms.}\label{subsection invariants de viterbo}

The invariants defined in this section have been introduced by C. Viterbo in \cite{V1}. We recall their construction. We first define invariants for Lagrangian submanifolds.

Let $L$ be an element of $\mathcal{L}$ and $S:M\times\R^q\to\R$ be one of its g.f.q.i's. Let us denote $S^{\lambda}=\{x\in
M\times\R^q\,|\,S(x)\leq\lambda\}$. Since $S$ is quadratic at infinity, the homotopy types of the pairs $(S^{\lambda},S^{\mu})$ and $(S^{\mu},S^{-\lambda})$ do not depend on $\lambda$, provided that $\lambda$ is sufficiently large . Therefore, we will denote $S^{\infty}$ and $S^{-\infty}$, instead of $S^{\lambda}$ and $S^{-\lambda}$ for $\lambda$ large enough.

Let us introduce $E^-_{\infty}$ the negative (trivial) bundle of the quadratic form which coincides with $S$ at infinity. We denote $B(E^-_{\infty})$, $S(E^-_{\infty})$ the ball bundle and the sphere bundle associated to $E^-_{\infty}$. The Thom isomorphism is given by $H^{\ast}(M)\to H^{\ast}(B(E^-_{\infty}),S(E^-_{\infty}))$, and we also have the isomorphism $H^{\ast}(B(E^-_{\infty}),S(E^-_{\infty}))\simeq H^{\ast}(S^{\infty},S^{-\infty})$. We will denote by $T$ their composition. For further informations on those isomorphisms, see \cite{H} for example. The inclusion $j_{\lambda}:S^{\lambda}\to S^{\infty}$ induces a morphism in cohomology $j_{\lambda}^{\ast}:H^{\ast}(S^{\infty},S^{-\infty})\to H^{\ast}(S^{\lambda},S^{-\infty})$, for all real number $\lambda$. We are now ready for the following.

\begin{definit}\label{def c}
Let $(u,L)\in H^{\ast}(M)\times\mathcal{L}$, with $u\neq 0$. Using a g.f.q.i $S$ of $L$, we define a real number $c(u,L)$ as follows:
\begin{equation}\label{eq c}
c(u,L)=inf\{\lambda\,|\,j_{\lambda}^{\ast}\circ T(u)\in H^{\ast}(S^{\lambda},S^{-\infty}) \text{ is non zero}\}.
\end{equation}
\end{definit}

Observe that $c(u,L)$ is well defined, and is independent of the
choice of $S$'s choice, up to additive constant. Indeed, if we
replace $S$ with an equivalent or stabilized generating function,
the value of $c(u,L)$ does not change, up to additive constant and
we conclude using theorem \ref{uniqueness}. Even if it doesn't
depend on the generating function, we sometimes use the notation
$c(u,S)$ instead of $c(u,L)$.

Since the cohomology of the sets $S^{\lambda}$ changes when we cross the level $c(u,L)$, it has to be a critical value of $S$.

Finally, observe that the definition can be extended to classes with compact support $u\in H_c^{\ast}(M)$.

\medskip
Then, we can use those invariants associated to Lagrangian submanifold to define other invariants associated to Hamiltonian diffeomorphisms.

Consider a Hamiltonian diffeomorphism $\psi\in\mathcal{H}(\RN)$. Its
graph $\Gamma_{\psi}$ is a Lagrangian submanifold of
$\overline{\R^{2n}}\times\R^{2n}$
($=(\R^{2n}\times\R^{2n},-\omega_0\oplus\omega_0)$, where $\omega_0$
is the standard symplectic structure on $\R^{2n}$). It coincides
with the diagonal $\Delta=\{(x,x)\,|\,x\in\R^{2n}\}$, outside the
product $B^{2n}(r)\times B^{2n}(r)$, for $r$ sufficiently large.
When we identify $\overline{\R^{2n}}\times\R^{2n}$ with
$T^{\ast}\Delta$ using the map,
$$(q,p,Q,P)\mapsto\left(\frac{q+Q}{2},\frac{p+P}{2},P-p,Q-q\right),$$
we see that the image $\widetilde{\Gamma_{\psi}}$ of $\Gamma_{\psi}$ is identified with the zero section of $T^{\ast}\Delta$ outside a compact set.

Then, we can associate the previous invariant to
$\widetilde{\Gamma_{\psi}}$ (We normalize generating functions by
asking their critical value at infinity to equal $0$). Let $1$ be a
generator of $H^0(\R^{2n})$ and $\mu$ a generator of
$H_c^{2n}(\R^{2n})$.

\begin{definit}[Viterbo, \cite{V1}]\label{def gamma}
We define,
$$c_-(\psi)=-c(\mu,\widetilde{\Gamma_{\psi}}),$$
$$c_+(\psi)=-c(1,\widetilde{\Gamma_{\psi}}),$$
$$\gamma(\psi)=c_+(\psi)-c_-(\psi),$$
$$\gamma(\phi,\psi)=\gamma(\psi^{-1}\phi).$$
\end{definit}

Let us describe now the properties of the numbers $\gamma$, $c_+$
and $c_-$ that we will use in the paper.

\begin{prop}[Viterbo, \cite{V1}]\label{propriétés gamma}

a)(Sign and Separation) For all $\psi$ in $\mathcal{H}$, we have
$$c_-(\psi)\leq 0\leq c_+(\psi).$$
Moreover, $c_-(\psi)=c_+(\psi)=0$ if and only if $\psi=\mathrm{Id}$.

\medskip
 b) (Triangle inequality) If $\phi$ is another diffeomorphism
in $\mathcal{H}$, then
$$c_+(\phi\circ\psi)\leq c_+(\phi)+c_+(\psi),$$
$$c_-(\phi\circ\psi)\geq c_-(\phi)+c_-(\psi),$$
$$\gamma(\phi\circ\psi)\leq \gamma(\phi)+\gamma(\psi).$$
$(\phi,\psi)\to\gamma(\phi,\psi)$ is a distance on $\mathcal{H}$.

\medskip
In particular, the separation property and the triangle inequality
imply that $(\phi,\psi)\to\gamma(\phi,\psi)$ is a distance on
$\mathcal{H}$.

\medskip
c)(Monotony) Let $\psi_1$ and $\psi_2$ be two Hamiltonians generated
by $H_1$ and $H_2$. Suppose that for all $(t,x)\in\R\times\R^{2n}$,
we have $H_1(t,x)\leq H_2(t,x)$. Then, $c_+(\psi_1)\leq c_+(\psi_2)$
and $c_-(\psi_1)\leq c_-(\psi_2)$.

As a consequence, if $H$ is a non-negative Hamiltonian, then
$c_-(\phi_H)=0$. If $H$ is in addition non zero, we deduce
$c_+(\phi_H)>0$.

\medskip
d) (Continuity) Let $H_1$ and $H_2$ be two compactly supported
hamiltonians, generating $\psi_1$ and $\psi_2$. Let $\|\cdot\|$ be
the usual norm on $C^0(\R^{2n}\times [0,1],\R)$. If $\Vert
H_1-H_2\Vert\leq\eps$, then
$|\gamma(\psi_1)-\gamma(\psi_2)|\leq\eps$.
\end{prop}

\subsection{Two symplectic capacities on $\R^{2n}$.}\label{def cap}

We start this section by reminding the reader of the definition of a
symplectic capacity. This is a "symplectic" way of measuring sets
 that plays an important role in symplectic topology. We will use it
 in particular for our convergence criterion in section
 \ref{critere}.

\begin{definit}[Ekeland-Hofer]\label{definition capacite}
A \textnormal{symplectic capacity} on $(\R^{2n},\omega_0)$ is a map associating to each subset $U\subset\R^{2n}$ a number $c(U)\in[0,+\infty]$ satisfying
\begin{enumerate}
\item $U\subset V\Rightarrow c(U)\leq c(V)$ (monotony),
\item $c(\phi(U))=c(U)$ for all Hamiltonian diffeomorphism $\phi\in\mathcal{H}$ (symplectic invariance),
\item $c(\lambda U)= \lambda^2c(U)$ for all real $\lambda>0$ (homogeneity),
\item $c(B^{2n}(1))=c(B^2\times\R^{2(n-1)})=\pi$, where $B^{2n}(1)$ is the unit ball of $\R^{2n}$ (normalisation).
\end{enumerate}
\end{definit}

The invariants defined in the previous section allow to define two symplectic capacities as follows (\cite{V1}).

\begin{definit} \textnormal{1}. For any compact subset $K\subset\R^{2n}$, we denote by $\gamma(K)$ the number defined by
$$\gamma(K)=\inf\{\gamma(\phi)\,|\,\phi(K)\cap K=\emptyset\}.$$
If $V$ is not compact, we set
$$\gamma(V)=\sup\{\gamma(K)\,|\,K\subset V\}.$$

\textnormal{2}. For any open subset $U\subset\R^{2n}$, we denote by $c(U)$ the number defined by
$$c(U)=\sup\{c_+(\phi_H)\,|\,\mathrm{Supp}(H)\subset U\}.$$
If $V$ is not an open set, we set
$$c(V)=\inf\{c(U)\,|\,V\subset U\}.$$
\end{definit}

The maps $c$ and $\gamma$ are symplectic capacities and moreover $c\leq\gamma$. We remind the reader of the definition of the displacement energy $$d(U)=\inf\{d_H(\phi,Id)\,|\,\phi(U)\cap U=\emptyset\},$$ where $d_H$ is Hofer's distance defined by
$$d_H(\phi,\psi)=\inf\{\|H-K\|\,|\,\text{$H$ generates $\phi$ and $K$ generates $\psi$}\},$$
whith $\|H\|=\int_0^1(\max_xH(t,x)-\min_xH(t,x))dt$.

We are going to define a new symplectic capacity derived from $c$,
but before we need the following lemma.

\begin{lem}\label{lemme capacité}
We consider a subset $V\subset\RN$ and $\R^2\times V\subset\R^{2+2n}$. Then, $$c(\R^2\times V)\geq c(V).$$
\end{lem}

That lemma follows from the reduction inequality of Proposition
\ref{red1}. We postpone its proof to Appendix. The reverse
inequality might be true but we are unable to prove it. That leads
us to introduce the following object.

\begin{definit} For any open subset $U\subset\R^{2n}$, we set
$$c^{\infty}(U)=\lim_{N\to\infty}c(U\times\R^{2N}),$$
and if $V$ is not an open subset,
$$c^{\infty}(V)=\inf\{c^{\infty}(U)\,|\,V\subset U\}.$$
\end{definit}

We obtain a symplectic capacity that satisfies
$c^{\infty}(V)=c^{\infty}(V\times\R^2)$ for all subset $V$ (this
property will be useful), and $c\leq c^{\infty}$. Moreover, since
$d(U)=d(U\times\R^{2k})$ and $c\leq d$, we have $c^{\infty}\leq d$.
To summarize the known inequalities between capacities we have,

\begin{prop}\label{capacites} $c\leq\gamma\leq d$ and $c\leq c^{\infty}\leq d$. \end{prop}

\subsection{Other distances derived from $\gamma$.}\label{distances derivées}

\medskip
In this section we introduce several other distances for many
reason. One is that we want to consider distances not only on
$\mathcal{H}$ but also on $Ham$. Another motivation is our result on
the Hamilton-Jacobi equation (section \ref{Hamilton Jacobi}) that
needs almost all of them. Finally, a stupid but important reason is
that we still don't know which is the best one to develop our
theory!

Let us start with the following distance defined on $\mathcal{H}$
already introduced by Cardin and Viterbo in \cite{CV}.

\begin{definit}\label{gamma tilde} For all Hamiltonian diffeomorphisms
$\phi,\psi\in\mathcal{H}$, we define
$$\tilde{\gamma}(\phi,\psi)=\sup\{\gamma(\psi^{-1}\phi(L)-L)\,|\,L\in\mathcal{L}\},$$
where $\gamma(L)=c(\mu,L)-c(1,L)$, $\forall L\!\in\!\mathcal{L}$ and
$L_1\!-\!L_2\!=\!\{(q,p_1-p_2)\,|\,(q,p_1)\in L_1,(q,p_2)\in L_2\}$,
for $L_1,L_2\in\mathcal{L}$.
\end{definit}

Then, we define distances not anymore on $\mathcal{H}$, but on
$Ham$.

\begin{definit} For any $H,K\in Ham$, we set $$\gamma_u(H,K)=\sup\{\gamma(\phi^t_{H},\phi^t_K)\,|\,t\in
[0,1]\}$$ and
$$\tilde{\gamma}_u(H,K)=\sup\{\tilde{\gamma}(\phi^t_{H},\phi^t_K)\,|\,t\in
[0,1]\}.$$
\end{definit}
Here, the subscript ``u'' means ``uniform''. Clearly, $\gamma_u$ and
$\tilde{\gamma}_u$ are distances on $Ham$.

For the next two distances, the principle is to add two dimensions
by associating to an Hamiltonian $H$ two suspensions defined on
$\R\times\R^{2+2n}$:
$$\hat{H}(s;t,\tau,x)=\tau+H(t,x),$$
$$\check{H}(s;t,\tau,x)=tH(st,x).$$
Here, the new time variable is $s$, while the former time variable
$t$ becomes a space variable (As a consequence $\hat{H}$ is an
autonomous Hamiltonian). We would like to define our distances by
$\hat{\gamma}(H,K)=\gamma(\hat{H},\hat{K})$ and
$\check{\gamma}(H,K)=\gamma(\check{H},\check{K})$. But since
$\hat{H}$ and $\check{H}$ are not compactly supported we have to be
slightly more subtle.

\begin{definit}\label{def suspendu} Let $\rho$ be a fixed real function defined on $[0,+\infty)$,  supposed to be non-negative, smooth, decreasing, with support in [0,1], flat at $0$ and such that $\rho(0)=1$. For every natural integer $\alpha$ and every real number $t$, we set $\rho_{\alpha}(t)=1$ if $-\alpha\leq t\leq\alpha$, and $\rho_{\alpha}(t)=\rho(|t|-\alpha)$ otherwise.

We denote by $\hat{H_{\alpha}}$ and $\check{H_{\alpha}}$ the
Hamiltonian functions defined on $\R\times\R^{2+2n}$, by
$$\hat{H_{\alpha}}(s;t,\tau,x)=\rho_{\alpha}(\tau)\tau+\rho_{\alpha}(t)H(t,x),$$
and $$\check{H_{\alpha}}(s;t,\tau,x)=\rho_{\alpha}(t)tH(st,x).$$

Then, for $H,K\in Ham$, we set
$$\hat{\gamma}(H,K)=\limsup_{\alpha\to+\infty} \gamma_u(\hat{H_{\alpha}},\hat{K_{\alpha}}),$$
and
$$\check{\gamma}(H,K)=\limsup_{\alpha\to+\infty} \gamma(\phi_{\check{H_{\alpha}}},\phi_{\check{K_{\alpha}}}).$$
\end{definit}

Remark that $\hat{\gamma}(H,K)$ and $\check{\gamma}(H,K)$ are
finite. Indeed, if we denote by $B$ a ball containing both supports
of $H$ and $K$, then $\hat{H}_{\alpha}$, $\hat{K}_{\alpha}$,
$\check{H}_{\alpha}$ and $\check{K}_{\alpha}$ have support in
$\R^2\times B$, for any integer $\alpha$. Hence
$\gamma(\phi_{\check{H_{\alpha}}},\phi_{\check{K_{\alpha}}})\leq 2
c(\R^2\times B)\leq 2c^{\infty}(B)$ (See section \ref{def cap} for
notations). It shows that the $\limsup$ in the definition of
$\check{\gamma}$ is finite. The same proof shows that
$\hat{\gamma}(H,K)$ is also finite.

The triangle inequality for $\hat{\gamma}$ and $\check{\gamma}$ is a
direct consequence of the triangle inequality for $\gamma$. The
separation property is obtained from the separation property for
$\gamma$ and Proposition \ref{inegalites distances}.

\medskip
For convenience, we will not write the subscript $\alpha$ anymore.
In the following, we will denote $\hat{H}$ for $\hat{H_{\alpha}}$,
and $\check{H}$ for $\check{H_{\alpha}}$.

\medskip
\noindent \textit{Remarks. } By repeating these constructions
several times (i.e., by taking suspensions of suspensions), we can
construct new distances. For example, we will use in section
\ref{Hamilton Jacobi} the distance
${{\gamma}}_2=\limsup_{\alpha\to+\infty}\check{\gamma}(\hat{H_{\alpha}},\hat{K_{\alpha}})$.

Using the invariance of $\gamma$, it is easy to verify that the
suspended distances $\hat{\gamma}$, $\check{\gamma}$ and $\gamma_2$
are invariant under the action of $\mathcal{H}$. Namely, for $H$,
$K$ Hamiltonians and $\varphi$ Hamiltonian diffeomorphism, we have:
$$\hat{\gamma}(H\circ\varphi,K\circ\varphi)=\hat{\gamma}(H,K),$$
$$\check{\gamma}(H\circ\varphi,K\circ\varphi)=\check{\gamma}(H,K),$$
$$\gamma_2(H\circ\varphi,K\circ\varphi)=\gamma_2(H,K).$$

\medskip
The following proposition gives inequalities between the distances.
It will be proved in Appendix. It is based on the reduction
inequality (Proposition \ref{red1}).
\begin{prop}\label{inegalites distances}
\begin{eqnarray*}
 & \tilde{\gamma}\leq\gamma,\\
 & \tilde{\gamma}_u\leq\gamma_u\leq\min(\hat{\gamma},\check{\gamma}).
\end{eqnarray*}
\end{prop}

\section{The convergence criterion.}\label{critere}

This is the central section of our paper. We give there the proof of
our main result, Theorem \ref{th intro}.

\subsection{A sufficient condition for a Hamiltonian diffeomorphism to be $\gamma$-close to $Id$.}

We start this section with some formulas concerning Hamiltonian flows. They
can be obtained by direct computation (see \cite{HZ}, page 144).

\begin{lem}\label{composition flots} For all Hamiltonians $H$ and $K$,
  with compact support, we have:
\begin{eqnarray*}
\phi_{\overline{H}}^t=(\phi_H^t)^{-1}, & \mathrm{where}\ \overline{H}(t,x)=-H(t,\phi^t(x))\\
\phi_{H\sharp K}^t=\phi_H^t\circ\phi_K^t,& \mathrm{where}\ (H\sharp K)(t,x)=H(t,x)+K(t,(\phi^t)^{-1}(x))\\
\phi_{\overline{H}\sharp K}^t=(\phi_H^t)^{-1}\circ\phi_K^t. &
\end{eqnarray*}
\end{lem}
\noindent
\textit{Remark.} $(\overline{H}\sharp K)(t,x)=(K-H)(t,\phi^t(x))$.

\medskip The following proposition shows that if a sequence of
Hamiltonians $(H_n)$ converges to zero uniformly on every compact
set contained in the complement of a set whose capacity is zero,
then $(\phi_{H_n})$ converges to $Id$ for $\gamma$.

\begin{prop}\label{premier lemme}
Let $H$ be a Hamiltonian on $\R^{2n}$ with compact support. If $U$ is an open subset of $\R^{2n}$, such that $c(U)\leq\eps$ and $|H(t,x)|\leq\eps$ for all $t\in [0,1]$ and all $x\in\R^{2n}-U$, then $\gamma(\phi_H)\leq 4\eps$.
\end{prop}
\textit{Proof.} Let $K_1$, $K_2$ be Hamiltonians with compact
support, such that $0\leq K_i\leq 1,\ i=1,2$, $K_1$ equals $1$ on
the support of $H$ and $K_2$ equals $1$ on the support of $K_1$
(hence $K_1\leq K_2$). Denote $\psi_{1,\eps}$ the diffeomorphism
generated by $H-\eps K_1$, and $\psi_{2,\eps}$ the diffeomorphism
generated by $\eps K_2$. Then we have $H\leq \eps K_2 + (H-\eps
K_1)$. As $(\psi_{2,\eps})^{-1}$ coincides with $Id$ on the support
of $H-\eps K_1$, the lemma \ref{composition flots} implies that
$\eps K_2 + H-\eps K_1$ is the Hamiltonian that generates
$\psi_{2,\eps}\circ\psi_{1,\eps}$. The monotony, the triangle
inequality and the continuity given Proposition \ref{propriétés
gamma} then give
$$c_+(\phi)\leq c_+(\psi_{2,\eps}\circ\psi_{1,\eps}) \leq c_+(\psi_{2,\eps})+ c_+(\psi_{1,\eps})\leq \eps + c_+(\psi_{1,\eps}).$$
Denote by $\widetilde{\psi_{1,\eps}}$ the diffeomorphism generated
by a non-negative Hamiltonian, with support in $U$, and greater than
$H-\eps K_1$. Then by the monotony property, $c_+(\psi_{1,\eps})\leq
c_+(\widetilde{\psi_{1,\eps}})$. Finally, since
$Supp(\widetilde{\psi_{1,\eps}})\subset U$, we get
$c_+(\widetilde{\psi_{1,\eps}})\leq c(U)\leq\eps$.

Using the inequality $H\geq -\eps K_2 + (H+\eps K_1)$, we obtain the same type of inequality for $c_-$.$\quad\Box$

\medskip
For example, if $K$ is a compact submanifold of dimension lower or equal than $n-1$, then $d(K)=0$ (and hence $c(K)=0$).

\subsection{What about non-identity elements that are close for $\gamma$?}

Unfortunately, the previous result cannot be straightforwardly generalised to obtain a general convergence criterions when the limit is not zero. Indeed, we can find two Hamiltonians that are $C^0$-close out of a null-capacity set, but not $\gamma$-close.

\medskip
\noindent \textit{Example. } It is well known that the capacities
$c$ and $\gamma$ of the unit sphere
$\mathcal{S}=\{x\in\R^{2n}\,|\,\|x\|=1\}$ are $\pi$. It is also true
for $c^{\infty}$. Then, for all $\eps>0$, there exists a Hamiltonian
$H$ with support in a small neighbourhood $U$ of $\mathcal{S}$, and
such that $c_+(\phi_H)>\pi-\eps$. Because of the monotony property
(proposition \ref{propriétés gamma}), $H$ can be chosen
non-negative. We set $U^+$ a neighbourhood of
$\{x\in\mathcal{S}\,|\,x_1\geq 0\}$ and $U^-$ a neighbourhood of
$\{x\in\mathcal{S}\,|\,x_1<0\}$, such that $U=U^+\cup U^-$. If $U$,
$U^+$ and $U^-$ are choosen small enough, we have $d(U^{\pm})<\eps$
and by proposition \ref{capacites} $c^{\infty}(U^{\pm})<\eps$. Using
some partition of unity associated to the decomposition $U=U^+\cup
U^-$, we get two functions $H^{\pm}$, with support in $U^{\pm}$ and
such that $H=H^++H^-$.

Now, we see that $H^+$ coincides with $H$ outside $U^-$, whose capacity verifies $c^{\infty}(U^{\pm})<\eps$, but on the other hand, $$\|\gamma(\phi_H)-\gamma(\phi_{H^+})\|\geq\gamma(\phi_H,\phi_{H^+})\geq\pi-\eps-\gamma(\phi_{H^-})\geq\pi-2\eps.$$
It shows that the previous statement is false when the limit is not zero.

\medskip
Nevertheless, we can introduce a new invariant, in order to extend the result of proposition \ref{premier lemme}.

\begin{definit} For any subset $U$ and any Hamiltonian $H\in Ham$, we define
$$\xi^{H}(U)=c^{\infty}\left(\bigcup_{t\in [0,1]}{\phi_H^t(U)}\right).$$
We may then set
$$\xi_{\lambda}(U)=\sup\,\xi^{H}(U), \text{
  for $0<\lambda\leq\infty$},$$
where the supremum is over all Hamiltonian functions $H$ with $\gamma_u(H)\leq\lambda$.
\end{definit}

\begin{theo}\label{le critere}
Let $H_1$ and $H_2$ be Hamiltonians on $\R^{2n}$ with compact support. Let $U$ be
a subset of $\R^{2n}$, satisfying one of the two following conditions:
\begin{enumerate}
\item $\xi_{\infty}(U)\leq\eps$.
\item $\exists\lambda>0$, $\xi_{\lambda}(U)=0$
\end{enumerate}
If $|H_1(t,x)-H_2(t,x)|\leq\eps$ for all $t\in[0,1]$ and all $x\notin
U$, then we have
$\gamma(\phi_{H_1},\phi_{H_2})\leq 4\eps$.
\end{theo}
\noindent
\textit{Proof.} Consider the Hamiltonian
$H(t,x)=H_1(t,\phi_2^t(x))-H_2(t,\phi_2^t(x))$. By assumption, $|H(t,x)|\leq\eps$, for all $(t,x)$ with
$x\notin\phi_2^{-t}(U)$ and hence for all $t$ and all
$x\notin\bigcup_{t\in [0,1]}{(\phi_2^{-1})^t(U)}$. Each
condition on $U$ implies $$c\left(\bigcup_{t\in
  [0,1]}{(\phi_2^{-1})^t(U)}\right)\leq c^{\infty}\left(\bigcup_{t\in
  [0,1]}{(\phi_2^{-1})^t(U)}\right)\leq\eps.$$
 By proposition \ref{premier
  lemme} and lemma \ref{composition flots}, we get $\gamma(\phi_{H_1},\phi_{H_2})=\gamma(\phi_H)\leq 4\eps$.$\quad\Box$

\medskip
\noindent \textit{Important remark. }In the proof of theorem \ref{le
critere}, we see that the important condition is in fact
$\xi^{H_2}(U)\leq \eps$, which is of course implied by both
conditions $\xi_{\infty}(U)\leq\eps$ and $\xi_{\lambda}(U)=0$.

\begin{corol}\label{critere pour les autres} The conclusion of theorem \ref{le critere} still holds if we replace $\gamma$ with $\tilde{\gamma}$.
For the distances on $Ham$, we get under the same assumptions
$$\mathfrak{d}(H_1,H_2)\leq 4 \|H_1-H_2\|_{C^0},$$
where $\mathfrak{d}$ is either $\tilde{\gamma_u}$, $\gamma_u$,
$\hat{\gamma}$ or $\check{\gamma}$.
\end{corol}
\noindent \textit{Proof.} By proposition \ref{inegalites distances}
(inequality between distances), we just have to prove it for
$\hat{\gamma}$ and $\check{\gamma}$. Then remark that under the
hypothesis of theorem \ref{le critere}, we have
$|\hat{H_1}(s;t,\tau,x)-\hat{H_2}(s;t,\tau,x)|\leq\eps$  and
$|\check{H_1}(s;t,\tau,x)-\check{H_2}(s;t,\tau,x)|\leq\eps$ for all
integer $\alpha$, all $s\in[0,1]$, and all
$(t,\tau,x)\notin\R^2\times U$.

Unfortunately, even if $U$ satisfies one of the conditions of proposition \ref{le critere}, it is not in general the case for $\R^2\times U$. However, by the above remark, it is sufficient to show that for all real number $\delta>0$ and all integer $\alpha$ large enough, $\xi^{{\check{H_2}}}(\R^2\times U)\leq \eps+\delta$ and  $\xi^{{\hat{H_2}}}(\R^2\times U)\leq \eps+\delta$. By letting $\delta$ tend to zero and taking limsup with respect to $\alpha$, we obtain $\hat{\gamma}(H,K)\leq 4\eps$ and $\check{\gamma}(H,K)\leq 4\eps$ as required.

Let us denote $F$ for $\check{H_2}$ or $\hat{H_2}$.
The inequalities on $\xi^{F}$ come directly from the expression of $\phi_{\check{H_2}}$ and $\phi_{\hat{H_2}}$ (see computations in Appendix \ref{suspended}). Indeed, in both cases,
$$\bigcup_{s\in [0,1]}{(\phi_F^{-1})^s([-\alpha,\alpha]^2\times U)}\subset\R^2\times \bigcup_{s\in [0,1]}{(\psi^{-1})^s(U)},$$
where $\psi^s$ is a Hamiltonian isotopy that appears in last
coordinate when we compute $\phi_F$. Therefore, since
$\xi^{F}(\R^2\times U)=
\lim_{\alpha\to+\infty}\xi^{F}([-\alpha,\alpha]^2\times U)$, we get
for any $\delta>0$ and any $\alpha$ large enough:
\begin{eqnarray*}
\xi^{F}(\R^2\times U) & \leq & \delta+\, c^{\infty}\!\left(\!\R^2\!\times\!\!\bigcup_{s\in [0,1]}{(\psi^{-1})^s(U)}\!\right)\\
& = & \delta+\, c^{\infty}\!\left(\bigcup_{s\in [0,1]}{(\psi^{-1})^s(U)}\!\right)\!\leq\delta+\eps.
\end{eqnarray*}
That concludes the proof.$\quad\Box$

\begin{corol}\label{suite} Let $(H_k)$ be a sequence of Hamiltonians in $Ham$, whose supports are contained in a fixed compact set. Suppose there exist a Hamiltonian $H\in Ham$ and a compact set $K\in\RN$ with $\xi_{\infty}(K)=0$, such that $(H_k)$ converges uniformly to $H$ on every compact set of $\R\times(\RN -K)$. Then $(\phi_{H_k})$ converges to $\phi_H$ for $\tilde{\gamma}$, $\gamma$, and $(H_k)$ converges to $H$ for $\tilde{\gamma_u}$, $\gamma_u$, $\hat{\gamma}$, $\check{\gamma}$.
\end{corol}
\noindent
\textit{Proof.} For  $\tilde{\gamma}$, $\gamma$, it is a direct consequence of the remark that follows theorem \ref{le critere}. We just have to verify that for all $\eps>0$, there exists a small neighbourhood $U$ of $K$ such that $\xi^{H}(U)\leq\eps$. This is true because for every neighbourhood $V$ of $\bigcup_{t\in[0,1]}\phi_H(K)$, we can choose a neighbourhood $U$ of $K$ such that  $$\bigcup_{t\in[0,1]}\phi_H(U)\subset V.$$ Since $c^{\infty}(\bigcup_{t\in[0,1]}\phi_H(K))=0$ and $\bigcup_{t\in[0,1]}\phi_H(K)$ is compact, we can choose $V$ such that $c^{\infty}(V)\leq\eps$, and obtain $c^{\infty}(\bigcup_{t\in[0,1]}\phi_H(U))\leq\eps$ as required.

For $\hat{\gamma}$ and $\check{\gamma}$, we have to verify that for all $\eps>0$ and all $\delta>0$, there exists a small neighbourhood $U$ of $K$ such that for all $\alpha$ large enough $\xi^{\phi}(U)\leq\eps+\delta$, where $F$ is either $\hat{H}$ or $\check{H}$. The proof made above for $\phi_H$ shows that we can find $U$ such that $\xi^{f}(U)\leq\eps$, where $f$ generates the isotopy $\psi^s$ defined as in the proof of corollary \ref{critere pour les autres}. Therefore we have for all $\delta$ and all $\alpha$ large enough, $\xi^{F}(\R^2\times U)\leq\xi^{F}([-\alpha,\alpha]^2\times U)+\delta\leq\xi^{f}(U)+\delta\leq\eps+\delta$.

By proposition \ref{inegalites distances}, corollary \ref{suite} is also true for $\tilde{\gamma_u}$ and $\gamma_u$.
$\quad\Box$

\medskip
\noindent
\textit{Remark. } Similar proofs give that theorem \ref{le critere} and corollary \ref{suite} still hold for ${{\gamma}}_2$.

\subsection{Example of a non trivial $\xi$-small set.}\label{xi
petit}

\begin{prop}\label{dim n-2} Let $U$ be a closed submanifold of $\R^{2n}$ whose dimension $d$ verifies $d\leq n-2$. Then $\xi_{\infty}(U)=0$.
\end{prop}

\textit{Proof. }Let $H\in Ham$. The problem is that $\bigcup_{t\in [0,1]}{\phi_H^t(U)}$ is not in general a manifold. To avoid that problem, we are going to add two dimensions and make a suspension in this way. We denote by $\Phi$ the Hamiltonian diffeomorphism on $\R^{2+2n}=\{(t,\tau,x)\}$ generated by the Hamiltonian $$[0,1]\times\R^{2+2n}\to\R,\ (s;t,\tau,x)\mapsto tH(ts,x).$$ We also set $V=\Phi([0,1]\times[-1,1]\times U)$. The computation of $\Phi$ gives
$$\Phi(t,\tau,x)=(t,\tau-H(t,x),\phi^t(x)).$$
We see that $\bigcup_{t\in [0,1]}{\phi_H^t(U)}$ can be obtained from $V$ by symplectic reduction by the coisotropic manifold $\{\tau=0\}$. So we are going to look for a Hamiltonian diffeomorphism $\phi_K$ that displaces $V$ and preserves $\{\tau=0\}$ at the same time. If the Hamiltonian does not depend on $t$, the second condition is verified. Since $V$ is compact, it is sufficient for $K$ to verify $$\forall v\in V,\ \R X_K(v)\cap T_vV=\{0\},$$ which is equivalent to $$\forall v\in V,\ \ker dK(v)\oplus
T_vV^{\perp}=\R^{2+2n}$$ and to $$\forall v\in V,T_vV^{\perp}\not\subset\ker dK(v).$$ That makes us consider the 1-jet bundle $J^1(\R\times\R^{1+2n},\R)$ and its submanifold
$$W=\{(s,q;\sigma,p;z)\,|\, (s,q)\in V, z\in\R, T_{(s,q)}V^{\perp}\subset \ker(\sigma,p)\}.$$
The dimension of $W$ is exactly $2n+1$. Indeed, the vector space $\{(\sigma,p)\in {\R^{2n+2}}^*\,|\, T_{(s,q)}V^{\perp}\subset \ker(\sigma,p)\}$ has dimension $2n+2- \text{dim}(T_{(s,q)}V^{\perp})= n$.

By Thom transversality theorem (see \cite{GG} for example), there exists a function $L$ whose 1-jet verifies $j^1L\pitchfork W$. But $j^1L$ can be seen as a function $\R\times\R^{1+2n}\to J^1(\R\times\R^{1+2n},\R)$, and by lemma 4.6 page 53 in \cite{GG}, we have for a generic choice of $s\in\R$,  $j^1L(s,\cdot)\pitchfork W$. We fix $s$ as previously and we denote $K:\R^{2+2n}\to\R$, $K(t,\cdot)=L(s,\cdot)$

Then, notice that for every $s,q,p,z$, the set of all $\sigma$ such that $(s,q;\sigma,p;z)\in W$ is either $\emptyset$ or $\R$. It can be shown by direct computation of $TV^{\perp}$, whose first component appears to be always $\{0\}$. As a consequence, we get $j^1K\pitchfork W$ ($j^1K$ differs from $j^1L(s,\cdot)$ just by its $\sigma$-component which is $\{0\}$ instead of $\derive{L}{s}(s,\cdot)$ for $j^1L(s,\cdot)$).

 Now, since $(2n+2)+(2n+1)=\text{dim}(j^1K(\R^{2+2n}))+\text{dim}(W)< \text{dim}(J^1(\R\times\R^{1+2n},\R))=4n+5$, we get $j^1K(\R^{2+2n})\cap W=\emptyset$. It follows that $K$ satisfies the two conditions: it preserves $\{\tau=0\}$ and it satisfies $$\forall v\in V,\ \R X_K(v)\cap T_vV=\{0\}.$$ As $V$ is compact, for $\eps$ small enough, since $\phi_{\eps K}=\phi_K^{\eps}$, we have $\phi_{\eps K}(V)\cap V=\emptyset$. In addition $\eps K$ can be made as $C^0$-small as we want.

We are now ready for the reduction by $\{\tau=0\}$. Since it preserves $\{\tau=0\}$, $\eps K$ induces a Hamiltonian on the reduction $\R^{2n}$. This Hamiltonian is $C^0$-small and generates a diffeomorphism $\psi$ whose Hofer's distance to identity $d_H(\psi,id)$ is small, and that satisfies $$\psi\left(\bigcup_{t\in [0,1]}{\phi_H^t(U)}\right)\cap \bigcup_{t\in [0,1]}{\phi_H^t(U)} =\emptyset .$$ This Hamiltonian is not compactly supported, but any Hamiltonian with compact support which coincides with it on a sufficiently large ball, would have the same properties. That proves $d\left(\bigcup_{t\in [0,1]}{\phi_H^t(U)}\right)=0$, and since $c^{\infty}\leq d$, we get $\xi^{H}(U)=0$. $\quad\Box$

\section{Completions and extension of Hamiltonian dynamics}\label{extension
dynamics}

In this section, we introduce the completions and give the first
properties of their elements: the existence of a flow that acts on
Lagrangian submanifolds, the notion of first integral and the
existence of a support. The full section \ref{Hamilton Jacobi} will
be devoted to another property related to the Hamilton-Jacobi
equation.

\subsection{Notations, inclusions and
definitions}\label{generalités}

Let us denote respectively $\overline{\mathcal{H}}^{\gamma}$,
$\overline{Ham}^{\gamma_u}$,
$\overline{\mathcal{H}}^{\tilde{\gamma}}$,
$\overline{Ham}^{\tilde{\gamma}_u}$,
$\overline{Ham}^{\hat{\gamma}}$, $\overline{Ham}^{\check{\gamma}}$
and $\overline{Ham}^{\gamma_2}$ the completions of
$(\mathcal{H},\gamma)$, $(Ham,\gamma_u)$,
$(\mathcal{H},\tilde{\gamma})$, $(Ham,\tilde{\gamma}_u)$,
$(Ham,\hat{\gamma})$, $(Ham,\check{\gamma})$ and $(Ham,\gamma_2)$.

The sets $\overline{\mathcal{H}}^{\gamma}$ and
$\overline{\mathcal{H}}^{\tilde{\gamma}}$ have a natural structure
of group with bi-invariant metric induced by the natural structures
on $(\mathcal{H},\gamma)$ and $(\mathcal{H},\tilde{\gamma})$.
Moreover we have the following fact:

\begin{prop}\label{laplication} The map $H\mapsto\phi_H^1$ induces Lipschitz maps
$\overline{Ham}^{\gamma_u}\to\overline{\mathcal{H}}^{\gamma}$ and
$\overline{Ham}^{\tilde{\gamma}_u}\to\overline{\mathcal{H}}^{\tilde{\gamma}}$.
\end{prop}

\noindent \textit{Proof. } Indeed, by construction of the distances,
$H\mapsto\phi_H^1$ is Lipschitz both as a map
$(Ham,\gamma_u)\to(\mathcal{H},\gamma)$ and as a map
$(Ham,\tilde{\gamma}_u)\to(\mathcal{H},\tilde{\gamma})$.$\quad\Box$

\medskip The inequalities between the different distances, proved in
Proposition \ref{inegalites distances}, induce inclusions between
the completions which may be summarized by the following diagram.
Here, $\overline{\mathcal{H}}^{d_H}$ denotes the completion of
$\mathcal{H}$ for Hofer's distance $d_H$ (which satisfies
$d_H\leq\gamma$) and $C_c$ the set of continuous (not necessarily
smooth) Hamiltonians with compact support.

$$\text{
\xymatrix{ C_c \ar@{^{(}->}[r] \ar@{^{(}->}[d] & \overline{Ham}^{\check{\gamma}} \ar@{^{(}->}[dr]\\
\overline{Ham}^{\gamma_2} \ar@{^{(}->}[r] & \overline{Ham}^{\hat{\gamma}} \ar@{^{(}->}[r] & \overline{Ham}^{\gamma_u} \ar@{^{(}->}[r] \ar[d] & \overline{Ham}^{\tilde{\gamma}_u} \ar[d] \\
 & \overline{\mathcal{H}}^{d_H} \ar@{^{(}->}[r] & \overline{\mathcal{H}}^{\gamma} \ar@{^{(}->}[r] & \overline{\mathcal{H}}^{\tilde{\gamma}}
}
}$$

As in Proposition \ref{laplication}, the map $(H,t)\mapsto\phi_H^t$,
$Ham\times\R\to\mathcal{H}$ induces maps
$\overline{Ham}^{\gamma_u}\times\R\to\overline{\mathcal{H}}^{\gamma}$
and
$\overline{Ham}^{\tilde{\gamma}_u}\times\R\to\overline{\mathcal{H}}^{\tilde{\gamma}}$.
Therefore, any element $H$ in one of the completions
$\overline{Ham}^{\gamma_u}$, $\overline{Ham}^{\tilde{\gamma}_u}$,
$\overline{Ham}^{\hat{\gamma}}$, $\overline{Ham}^{\check{\gamma}}$
or $\overline{Ham}^{\gamma_2}$ can be associated a path in either
$\overline{\mathcal{H}}^{\gamma}$, or
$\overline{\mathcal{H}}^{\tilde{\gamma}}$. This path of course has
the semi-group property. That leads us to the following definition.

\begin{definit}Such a path will be called the \textnormal{generalized Hamiltonian flow} generated by $H$.
\end{definit}

\subsection{Action on Lagrangian submanifolds}\label{section action
lag}

Recall that the set $\mathcal{L}$ of Lagrangian submanifolds
isotopic to the zero section by compactly supported Hamiltonian
isotopy, can be endowed with Viterbo's distance,  also denoted
$\gamma$ (set $\gamma(L_1,L_2)=\gamma(L_1-L_2)$, see definition
\ref{gamma tilde}). Let us denote $\mathfrak{L}$ the completion of
$\mathcal{L}$ with respect to this distance.

\begin{prop}\label{action lag} The groups $\overline{\mathcal{H}}^{\gamma}$ and $\overline{\mathcal{H}}^{\tilde{\gamma}}$
naturally act on the set $\mathfrak{L}$. This action extends the
natural action of $\mathcal{H}$ on $\mathcal{L}$.
\end{prop}

\medskip
\noindent \textit{Proof. } It is a simple consequence of the
inequality $\tilde{\gamma}\leq\gamma$ (Proposition \ref{inegalites
distances} proved in Appendix).

Let $L\in\mathfrak{L}$ represented by a sequence $(L_k)$ and $\phi$
in $\overline{\mathcal{H}}^{\gamma}$ (the proof is the same for
$\overline{\mathcal{H}}^{\tilde{\gamma}}$), represented by a
sequence $(\phi_k)$. We are going to show that $(\phi_k(L_k))$
defines an element of $\mathfrak{L}$ that we will denote $\phi(L)$.

This follows easily from the fact that for $\phi, \psi\in
\mathcal{H}$ and $L,M\in\mathcal{L}$,
\begin{eqnarray*}\gamma(\phi(L)-\psi(M)) & \leq & \gamma(\phi(L)-\psi(L)) +
\gamma(\psi(L)-\psi(K))\\
 & \leq  & \gamma(\psi^{-1}\phi(L)-L)  + \gamma(L-K)\\
 & \leq  & \tilde{\gamma}(\phi,\psi) + \gamma(L-K).\quad\Box
\end{eqnarray*}

\medskip
\noindent \textit{Remark. } A consequence of Proposition \ref{action
lag} is that we can define what is a Lagrangian submanifold
\textit{invariant} under a generalized flow.

That leads us to another question which is: Can we define what is an
invariant hypersurface of a generalized flow?

A (partial) answer to this question is that we can define what is a
first integral of a generalized Hamiltonian flow.

\subsection{Notion of first integral}

This property has been first mentioned in \cite{CV}, in the
definition (3.3) of the so-called $c$-commuting Hamiltonians. Let us
restate it with our notations.

An element in one of the completions $\overline{Ham}^{\gamma_u}$,
$\overline{Ham}^{\tilde{\gamma}_u}$,
$\overline{Ham}^{\hat{\gamma}}$, $\overline{Ham}^{\check{\gamma}}$
and $\overline{Ham}^{\gamma_2}$ will be said \textit{autonomous} if
it can be represented by a Cauchy sequence of time-independent
Hamiltonian functions.

\begin{definit} Let $H, K$ be two elements in one of the above
completions, generating two respective generalized flows $\phi_H^t$
and $\phi_K^t$. Then we will say that $H$ and $K$
\textnormal{commute}, or that $K$ is a \textnormal{first integral}
of $H$ if $\phi_K^s\phi_H^t\phi_K^{-s}\phi_H^{-t}=Id$.

In other words, $K$ is a \textnormal{first integral} of $H$ if there
exists two Cauchy sequences $(H_n)$ and $(K_n)$ representing $H$ and
$K$, such that for all $s$ and $t$,
$\phi_{K_n}^s\phi_{H_n}^t\phi_{K_n}^{-s}\phi_{H_n}^{-t}$ c-converges
to $Id$.
\end{definit}

It is proved in \cite{CV} that this definition extends the usual
definition of commuting Hamiltonian functions.

\subsection{Existence of a support}

In this section, we state a lemma which makes it possible to define a support for the elements of the different completions.

\begin{lem}\label{support}
\begin{description}
\item[a. ]Let $(\phi_n)$ be a sequence in $\mathcal{H}$ converging to a Hamiltonian diffeomorphism $\phi$, with respect to $\gamma$ or $\tilde{\gamma}$. Assume that there exists a set $U\in\RN$ such that $\text{supp}(\phi_n)\subset U$. Then $\text{supp}(\phi)\subset \overline{U}$.
\item[b. ]Let $(H_n)$ be a sequence in $Ham$ converging to a smooth Hamiltonian function $H$, with respect to $\gamma_u$, $\tilde{\gamma}_u$, $\hat{\gamma}$, $\check{\gamma}$, etc. Assume that there exists a set $U\in\RN$ such that $\text{supp}(H_n)\subset U$. Then $\text{supp}(H)\subset \overline{U}$.
\end{description}
\end{lem}
\textit{Proof. } \textbf{a. } Thanks to Proposition \ref{inegalites distances}, we just have to prove the assertion in the case of $\tilde{\gamma}$. Suppose $\text{supp}(\phi)\not\subset\overline{U}$. Then there exists an $x$ in $\RN-\overline{U}$ such that $\phi(x)\not=x$. Let $\psi$ be a Hamiltonian diffeomorphism whose support is included in $\RN-\overline{U}$ and which does not contain $\phi(x)$. Suppose in addition that $\psi(x)\not=x$. Then, since the supports of $\phi_n$ and $\psi$ are disjoint, we have $\psi\circ\phi_n^{-1}\circ\psi^{-1}\circ\phi_n=Id$, for all integer $n$. Taking limit, we get on one hand $\psi\circ\phi^{-1}\circ\psi^{-1}\circ\phi=Id$. But on the other hand, we have by construction, $\psi\circ\phi^{-1}\circ\psi^{-1}\circ\phi(x)=\psi(x)\not=x$, which is contradictory.

\noindent
\textbf{b. } We use the first part of the lemma to conclude that for all time $t$, $\text{supp}(\phi^t)\subset \overline{U}$. This implies that $\text{supp}(H)\subset \overline{U}$.$\quad\Box$

\medskip
\textit{Remark. }A similar argument shows that the property of letting globally invariant any sphere centered at $0$, is invariant by taking  $\gamma$ or $\tilde{\gamma}$ limits. Similarly, a $\gamma_u$, $\tilde{\gamma}_u$, $\hat{\gamma}$ or $\check{\gamma}$ limit of radial Hamiltonians is radial.

\begin{definit}\label{def support}
\textbf{a. }Let $\psi$ be an element of
$\overline{\mathcal{H}}^{\gamma}$ or
$\overline{\mathcal{H}}^{\tilde{\gamma}}$. Then we define
$\text{support}(\psi)$ as
$$\bigcap\{\overline{U}\,|\,\text{$U$ open set, such that}\text{ there exists $(\psi_n)$ representing $\psi$ such that }$$
$$\quad\quad\quad\quad\quad\quad\quad\quad\quad\quad\quad\quad\quad\quad\quad\quad\quad\quad\quad\forall n,\text{supp}(\psi_n)\subset U\},$$
where ``$\text{supp}$'' denotes the usual notions of support for
Hamiltonian diffeomorphisms.

\medskip \noindent
\textbf{b. }Let $K$ be an element of $\overline{Ham}^{\gamma_u}$,
$\overline{Ham}^{\tilde{\gamma}_u}$,
$\overline{Ham}^{\hat{\gamma}}$, $\overline{Ham}^{\check{\gamma}}$
or $\overline{Ham}^{\gamma_2}$. Then we define $\text{support}(K)$
as
$$\bigcap\{\overline{U}\,|\,\text{$U$ open set, such that}\text{ there exists $(K_n)$ representing $K$ such that }$$
$$\quad\quad\quad\quad\quad\quad\quad\quad\quad\quad\quad\quad\quad\quad\quad\quad\quad\quad\quad\forall n,\text{supp}(K_n)\subset U\},$$
where ``$\text{supp}$'' denotes the usual notions of support for
smooth Hamiltonians.

\end{definit}

These new notions of support coincide with the usual notions for
smooth Hamiltonians and Hamiltonian diffeomorphisms. Indeed, let
$\eta$ be either a Hamiltonian diffeomorphism viewed as an element
of $\overline{\mathcal{H}}^{\gamma}$ or
$\overline{\mathcal{H}}^{\tilde{\gamma}}$, or a smooth Hamiltonian
seen as an element of $\overline{Ham}^{\gamma_u}$,
$\overline{Ham}^{\tilde{\gamma}_u}$,
$\overline{Ham}^{\hat{\gamma}}$, $\overline{Ham}^{\check{\gamma}}$
or $\overline{Ham}^{\gamma_2}$. Let $(\eta_n)$ be a sequence
representing $\eta$, and $U$ an open set with
$\text{supp}(\eta_n)\subset U$ for all $n$. Then lemma \ref{support}
gives $\text{supp}(\eta)\subset\overline{U}$. Hence
$\text{supp}(\eta)\subset\text{support}(\eta)$. Conversely, for any
neighbourhood $U$ of $\text{supp}(\eta)$ the constant sequence
$(\eta)$ converges to $\eta$ and has support in $U$. Therefore
$\text{support}(\eta)\subset\bigcap_{\mathcal{V}}\overline{U}$,
where the intersection is over the set $\mathcal{V}$ of all open
neighbourhoods of $\text{supp}(\eta)$. Then, it is easy to see that
$\bigcap_{\mathcal{V}}\overline{U}=\bigcap_{\mathcal{V}}U=\text{supp}(\eta)$.

\section{Description of some elements of the
completions}\label{description elements}

The elements of the different completions are by definition
equivalence classes of Cauchy sequences. So they are defined in a
very abstract way. In this section, we show that many elements of
the completions can be seen in a more concrete way.

\subsection{Examples in the completion of $Ham$}

The inequalities between Hofer's distance and our four distances
$\gamma_u$, $\tilde{\gamma}_u$, $\hat{\gamma}$ and $\check{\gamma}$
on $Ham$ imply inclusions of the completions. In particular any
continuous time-dependent Hamiltonian can be seen as an element of
$\overline{Ham}^{\gamma_u}$, $\overline{Ham}^{\tilde{\gamma}_u}$,
$\overline{Ham}^{\hat{\gamma}}$ and
$\overline{Ham}^{\check{\gamma}}$.

In view of Theorem \ref{le critere} and Corollary \ref{suite}, we
can conjecture that if $K\in\RN$ satisfies $\xi^{\infty}(K)=0$ a
sequence $(H_k)$ converges uniformly on compact sets of
$\R\times(\RN-K)$ to a function $H$ continuous on $\R\times(\RN-K)$,
then $(H_k)$ is Cauchy for either $\gamma_u$, $\tilde{\gamma}_u$,
$\hat{\gamma}$ or $\check{\gamma}$ (compare with Corollary
\ref{suite}). We are still unable to prove it, but if we restrict to
a family of Hamiltonians which converge to $+\infty$ at their
discontinuity points, this result can be established.

\begin{definit}\label{la famille autonome}
We denote by $\mathfrak{F}$ the set of all functions
$H:\R\times\RN\to\R\cup\{+\infty\}$ such that:
\begin{description}
\item[(i)]  There exist $K\in\RN$ with $c^{\infty}(K)=0$ such that $H(t,x)=+\infty\Rightarrow x\in K$,
\item[(ii)] $H$ vanishes at infinity: $\forall\eps>0,\exists r,(|x|>r\Rightarrow (\forall t, |H(t,x)|<\eps))$,
\item[(iii)] $H$ is continuous on $\R\times\RN$.
\end{description}
We also set $\mathfrak{F}^{\infty}=\{H\in\mathfrak{F}\,|\,\text{$H$
is smooth on $\R\times\RN-H^{-1}(\{+\infty\})$}\}$, and
$\mathfrak{A}$, $\mathfrak{A}^{\infty}$ the subsets of
time-independent elements of  $\mathfrak{F}$ and
$\mathfrak{F}^{\infty}$.
\end{definit}

For the elements of $\mathfrak{F}^{\infty}$, the set of
discontinuity is somehow "stable" under the Hamiltonian flow. This
property allows to consider functions with a larger discontinuity
set than what could be expected in the general case
($c^{\infty}(K)=0$ instead of $\xi^{\infty}(K)=0$).

\begin{lem}\label{c dedans autonome} Suppose $H$ is an element of $\mathfrak{A}$ and $K=H^{-1}(\{+\infty\})$. Then there exists a sequence of smooth autonomous Hamiltonians $(H_k)\in Ham$ with the following properties:
\begin{description}
 \item[a.] $(H_k)$ converges to $H$ uniformly on every compact subset of $\RN-K$.
 \item[b.] $(H_k)$ is Cauchy for $\gamma_u$, $\tilde{\gamma}_u$, $\hat{\gamma}$ and $\check{\gamma}$.
\end{description}
Moreover, if $H\in\mathfrak{A}^{\infty}$, then any sequence $(H_k)$
that converges to $H$ uniformly on the compact subsets of $\RN-K$,
does not converge in $Ham$, for none of the distances $\gamma_u$,
$\tilde{\gamma}_u$, $\hat{\gamma}$ and $\check{\gamma}$.
\end{lem}

\medskip
\noindent \textit{Proof. } Fix $k>0$. Properties (ii) and (iii) in
Definition \ref{la famille autonome} imply that $K$ is compact.
Since $c^{\infty}(K)=0$, there exists an open neighborhood $U$ of
$K$ such that $c^{\infty}(U)\leq \frac1k$. Then, if we denote
$H^{>A}=\{x\,|\,H(x)>A\}$, we have for $A$ large enough, $K\subset
H^{>A}\subset U$. Indeed, if it was not true, then for all integer
for all integer $l>0$, there would exists a point $a_l$ in $H^{>l}$,
but not in $U$. Then, the sequence $(a_l)$ would take values in
$H^{\geq 1}\cap (\RN-U)$ which is compact, and hence it would have a
subsequence that would converge to an element of $K\cap(\RN-U)$,
which contradicts our assumption. Let us fix a real number $A_k$
such that $H^{>A_k}\subset U$.

Now, let $H_k$ be a smooth function with compact support such that
$|H_k-H|<\frac{1}{k}$ on $\RN-H^{>A_k+\frac{2}{k}}$, and such that
$|H_k-A_k-\frac{2}{k}|<\frac{1}{k}$ on $H^{>A_k+\frac{2}{k}}$. The
sequence $(H_k)$ clearly converges to $H$ uniformly on every compact
subset of $\RN-K$. Let us see why it is Cauchy.

By Proposition \ref{inegalites distances}, we just have to prove it
for $\hat{\gamma}$ and $\check{\gamma}$. We write $F_k$ for either
$\check{H_k}$ or $\hat{H_k}$. We also denote, as in the proof of
Corollary \ref{critere pour les autres}, $\psi_k$ for the third
coordinate of $\phi_{F_k}$. Since $H_k$ is an autonomous
Hamiltonian, its flow $\phi_{H_k}^t$ preserves its level sets.
Hence, the isotopy $\psi_k^s$ preserves the level sets of $F_k$ (see
the computations in Appendix \ref{suspended}). Therefore, since by
construction $H^{>A_k+\frac{2}{k}}\subset H_k^{>A_k+\frac{1}{k}}$,
we have
$$\bigcup_{t\in[0,1]}\psi_k^t(H^{>A_k+\frac{2}{k}})\subset H_k^{>A_k+\frac{1}{k}}.$$
Let $\delta>0$ and suppose $\alpha$ is sufficiently large. Then, as
in the proof of Corollary \ref{critere pour les autres},

\begin{eqnarray*}
\xi^{F_k}(\R^2\times H^{>A_k+\frac{2}{k}}) & \leq & \delta+\, c^{\infty}\!\left(\!\R^2\!\times\!\bigcup_{s\in [0,1]}{(\psi_k^{-1})^s(H^{>A_k+\frac{2}{k}})}\!\right) \\
 &  \leq  & \delta+\, c^{\infty}\!\left(\bigcup_{s\in [0,1]}{(\psi_k^{-1})^s(H^{>A_k+\frac{2}{k}})}\!\right)\\
 & \leq & \delta+ c^{\infty}(H_k^{>A_k+\frac{1}{k}}).
\end{eqnarray*}
Since $H_k^{>A_k+\frac{1}{k}}\subset H^{>A_k}\subset U$ and
$c^{\infty}(U)\leq\frac1k$, we obtain
$\xi^{F_k}(H^{>A_k+\frac{2}{k}})\leq\frac1k+\delta$.

Now, pick an integer $l\geq k$. If $l$ and $k$ are large enough,
then we have $|\hat{H}_k-\hat{H}_l|\leq\frac1k$ and
$|\check{H}_k-\check{H}_l|\leq\frac1k$ on $\R^{2+2n}-(\R^2\times
H^{>A_k+\frac{2}{k}})$. Therefore, by the remark that follows
Theorem \ref{le critere}, we get $\hat{\gamma}(H_l,H_k)\leq \frac4k$
and $\check{\gamma}(H_l,H_k)\leq \frac4k$, after taking limsup with
respect to $\alpha$. It proves that $(H_k)$ is a Cauchy sequence for
$\tilde{\gamma}_u$, $\gamma_u$, $\hat{\gamma}$, and
$\check{\gamma}$.

Suppose now that $H$ is smooth on $\RN-K$. Then we can choose $H_k$
such that it coincides with $H$ on $B_k-H^{>A_k+\frac{2}{k}}$, where
$B_k$ is the ball of radius $k$, centered at $0$. Suppose that
$(H_k)$ converges to a Hamiltonian $L\in Ham$ for
$\tilde{\gamma}_u$, $\gamma_u$, $\hat{\gamma}$, and
$\check{\gamma}$. Then for any integer $k$, $\overline{H_k}\sharp
H_l$ converges to $\overline{H_k}\sharp L$ while $l$ tends to
infinity for $\tilde{\gamma}_u$ (see Lemma \ref{composition flots}
for notations). According to Lemma \ref{support}, since
$\overline{H_K}\sharp H_l$ has support in the complementary of
$B_k-H^{>A_k+\frac{2}{k}}$, $\overline{H_k}\sharp L$ has support in
its closure and hence $L$ coincides with $H$ on
$B_k-H^{>A_k+\frac{2}{k}}$. Since it is true for any $k$, $L$ has to
coincide with $H$ on $\RN-K$. Therefore $L$ cannot belong to $Ham$,
which contradicts our assumptions.

Finally, if $(L_k)$ is another sequence of Hamiltonians that
converges to $H$ uniformly on the compact subsets of $\RN-K$, then,
similarly as in the above proof that $(H_k)$ is Cauchy, we obtain
that $\hat{\gamma}(L_k,H_k)$ and $\check{\gamma}(L_k,H_k)$ converge
to $0$, where $H_k$ is the particular sequence defined in the
previous paragraph. Since $(H_k)$ does not converge, $(L_k)$ does
not converge either.$\quad\Box$

\medskip
\noindent \textit{Remark. } As usual, the results of Lemma \ref{c
dedans autonome} still hold for ${\gamma}_2$.

\begin{prop}\label{inclusion autonome} The set $\mathfrak{F}^{\infty}$ can be embedded into each completion $\overline{Ham}^{\gamma_u}$, $\overline{Ham}^{\tilde{\gamma}_u}$, $\overline{Ham}^{\hat{\gamma}}$ and $\overline{Ham}^{\check{\gamma}}$.
\end{prop}

\noindent \textit{Proof. } Let us first consider the autonomous case
(elements of $\mathfrak{A}^{\infty}$).

Since
$\overline{Ham}^{\hat{\gamma}}\subset\overline{Ham}^{\gamma_u}\subset\overline{Ham}^{\tilde{\gamma}_u}$
and
$\overline{Ham}^{\check{\gamma}}\subset\overline{Ham}^{\gamma_u}\subset\overline{Ham}^{\tilde{\gamma}_u}$,
it is enough to prove it for $\hat{\gamma}$ and $\check{\gamma}$. We
will make the proof for $\hat{\gamma}$ and the proof for
$\check{\gamma}$ will be exactly the same.
 Let $J$ be the function that associates to any $H\in\mathfrak{A}^{\infty}$ the element of $\overline{Ham}^{\hat{\gamma}}$ represented by any sequence $(H_k)$ that converges uniformly to $H$ on the compact sets of $\RN-H^{-1}(\{+\infty\})$. As we noticed at the end of the proof of Lemma \ref{c dedans autonome}, two such sequences are equivalent and hence $J$ is well-defined.

Let us now prove that $J$ is one-one. Let
$H,G\in\mathfrak{A}^{\infty}$ and let $(H_k)$, $(G_k)$ be two
sequences respectively associated to them, precisely constructed as
in the last but one paragraph of the previous proof. Suppose that
$G\neq H$, we are going to show that $\gamma(H_k,G_k)$ does not
converge to zero, that will imply that $\hat{\gamma}(H_k,G_k)$ does
not converge to zero.

We can define almost everywhere the flows $\phi_G^t$, $\phi_H^t$ and
$\psi^t=\phi_{G}^{-t}\circ\phi_{H}^t$. Let
$\psi_k=\phi_{G_{k}}^{-1}\circ\phi_{H_{k}}$. Since $G\neq H$, there
exists a point $x$ such that $\psi(x)\neq x$. Hence, there exists a
small ball $B$ around $x$ such that $\psi(B)\cap B=\emptyset$. Let
$K$ be a compact neighborhood of $\bigcup_t\psi^t(B)$. For $k$ large
enough, $H_k$ and $G_k$ coincide respectively with $H$ and $G$ on
$K$, and thus $\psi_k(B)\cap B=\emptyset$ too. Since
$\gamma(H_k,G_k)=\gamma(\psi_k)\geq\gamma(B) >0$, $\gamma(H_k,G_k)$
cannot converge to zero.

To achieve the proof, we just have to notice that the map $H\mapsto
\hat{H}$ is a one-one map
$\mathfrak{F}^{\infty}_{2n}\to\mathfrak{A}^{\infty}_{2n+2}$ (the
subscript denotes the dimension of the ambient symplectic space).
Thus, according to the autonomous case, if $H$ is in
$\mathfrak{F}^{\infty}_{2n}$ then $\hat{H}$ is in
$\overline{Ham}^{\gamma_u}(\R^{2n+2})$. Moreover, according to Lemma
\ref{c dedans autonome}, we may construct a Cauchy sequence $H_k$ of
the form $\hat{F_k}$ for some Hamiltonians $F_k$. That means that
$H$ is actually an element of
$\overline{Ham}^{\hat{\gamma}}(\R^{2n})$. Inclusions between
completions give that it is an element of
$\overline{Ham}^{\gamma}(\R^{2n})$ and
$\overline{Ham}^{\tilde{\gamma}}(\R^{2n})$ too. finally, a similar
reasoning using $\gamma_2$ instead of $\hat{\gamma}$ allows to see
$H$ as an element of
$\overline{Ham}^{\check{\gamma}}(\R^{2n})$.$\quad\Box$

\medskip
Now, if we denote by $C_0$ the set of continuous Hamiltonians that
vanish at infinity, we can improve the diagram of section
\ref{generalités}:
$$\text{
\xymatrix{\mathfrak{F} \ar@/_/[dd] \ar@/^/[drr]\\
 & \mathfrak{F}^{\infty}\cup C_0 \ar@{^{(}->}[r] \ar@{^{(}->}[dl] & \overline{Ham}^{\check{\gamma}} \ar@{^{(}->}[d]\\
\overline{Ham}^{\gamma_2} \ar@{^{(}->}[r] & \overline{Ham}^{\hat{\gamma}} \ar@{^{(}->}[r] & \overline{Ham}^{\gamma_u} \ar@{^{(}->}[r] \ar[d] & \overline{Ham}^{\tilde{\gamma}_u} \ar[d] \\
 & \overline{\mathcal{H}}^{d_H} \ar@{^{(}->}[r] & \overline{\mathcal{H}}^{\gamma} \ar@{^{(}->}[r] & \overline{\mathcal{H}}^{\tilde{\gamma}}
} }$$

\subsection{Examples in the completions of $\mathcal{H}$}

In the completions of $Ham$ easy examples was given by continuous
Hamiltonian functions. In the completions of $\mathcal{H}$ there are
no similar result. Indeed, there are no known relation between
$C^0$-distance and $\gamma$.

However, we can give concrete examples of elements of the completion
of  $\mathcal{H}$ by Corollary \ref{inclusion autonome}. Indeed, it
implies that the (generalized) flows generated by elements of
$\mathfrak{F}^{\infty}$ are in both
$\overline{\mathcal{H}}^{\gamma}$ and
$\overline{\mathcal{H}}^{\tilde{\gamma}}$. Let us give some examples
(in their construction, $\gamma$ can be replaced by $\tilde{\gamma}$
without any problem).

\medskip
\noindent \textit{Example of a non smooth homeomorphism in
$\overline{\mathcal{H}}^{\gamma}$.}

\medskip
\noindent We consider a decreasing function
$h:[0,+\infty)\to[0,+\infty)$, with support in $[0,1]$, and equal to
$1$ on $[0,3/4]$. Then we define $H_k(x)=\sum_{i=1}^k{h(2^i|x|^2)}$,
for $x\in\RN$ and $H(x)=\sum_{i=1}^{\infty}{h(2^i|x|^2)}$ (the
"sky-scrapper" Hamiltonian). Let us see why $\phi_H$ can be seen as
a non-smooth homeomorphism.

Lemma \ref{support} implies that $\phi_H$ coincides with
$\phi_{H_k}$ out of $B_{2^{-k}}$. So, we can compute the explicit
form of $\phi$. In polar coordinates, we obtain:
$$\phi_H(r,\theta)=(\theta-r^2f'(r^2),r),$$
for $r>0$ where $f(s)=\sum_{i>0}h(2^is)$ (for any $s$, all the terms
in this sum are 0 except maybe one). We see that $\phi_H$ is a
homeomorphism. Let us prove that it cannot  be smooth at zero.

 If we denote by $(q,p)$ the
coordinates in $\RN$, and by $\phi_1$ the projection of $\phi_H$ on
$\R^n\times\{0\}$, we have for $q\in 2^{-i}[1/2,1]$,
\begin{equation}\label{equa}
\derive{\phi_1}{q}(q,0)=\cos(q^22^ih'(2^iq^2))-2(q^42^{2i}h''(2^iq^2)+q^32^ih'(2^iq^2))\sin(q^22^ih'(2^iq^2)).
\end{equation}
Suppose that $h$ is chosen so that there exists $q_1$ and $q_1'$ in
$[1/2,1]$ such that
$\derive{\phi_1}{q}(q_1,0)\not=\derive{\phi_1}{q}(q_1',0)$ (we
denote by $A$ their difference), and define $q_i=\sqrt{2^{-i}}q_1$
and $q_i'=\sqrt{2^{-i}}q_1'$. Then, $(q_i)$ and $(q_i')$ vanish, but
from (\ref{equa}) we see that
$\derive{\phi_1}{q}(q_i,0)-\derive{\phi_1}{q}(q_i',0)$ converges to
$A$. Therefore $\phi_H$ cannot be smooth at $0$.

\medskip
\noindent \textit{Example of a discontinuous element in
$\overline{\mathcal{H}}^{\gamma}$.}

\medskip
\noindent In the previous example, the sequence of diffeomorphisms
$(\phi_{H_k})$ was converging almost everywhere to a homeomorphism
(which was not a diffeomorphism). Therefore, one could think that
the class of $(\phi_{H_k})$ in the completion
$\overline{\mathcal{H}}^{\gamma}$ can be represented by a
homeomorphism. However, with the help of Proposition \ref{c dedans
autonome}, we can show that it is not true in general, at least in
dimension $2n\geq 4$.

Indeed, consider $H:\R^2\times\RN\to \R$, $$(x_1,x_2)\mapsto
\frac{1}{|\|x_1\|^2-1|+\|x_2\|^2}\,\chi(\|(x_1,x_2)\|),$$ where
$\chi$ is smooth with compact support and equals $1$ on the ball of
radius $2$ centered at zero. Clearly, $H\in\mathfrak{F}^{\infty}$
(because $K=H^{-1}(\{+\infty\})=\S^1\times\{0\}$ satisfies
$c^{\infty}(K)=0$ as required). Consider the sequence $(H_k)$
constructed in the proof of Lemma \ref{c dedans autonome}. Since
$(H_k)$ is Cauchy for $\gamma_u$, $(\phi_{H_k})$ is Cauchy for
$\gamma$. Suppose it converges to an element $\phi$. Then, Lemma
\ref{support} implies that for any neighbourhood $U$ of $K$ and for
$k$ large enough, $\phi$ coincides with $\phi_{H_k}$ on
$\R^{2+2n}-U$. Therefore, we can compute the explicit form of $\phi$
on $\R^{2+2n}-K$.

In polar coordinates $(s_1,\theta_1,s_2,\theta_2)$ with
$s_1=\|x_1\|^2$ and $s_2=\|x_2\|^2$, we get for $s_1<1$:
$$\phi(s_1,0,0,0)=\left(s_1,\frac{s_1}{(1-s_1)^2},0,0\right).$$
If we let $s_1$ converge to $1$, we see that $\phi$ is not
continuous.

\medskip
\noindent \textit{Questions. } The previous examples lead us to
natural questions: Are all the elements of
$\overline{\mathcal{H}}^{\gamma}(\R^2)$ homeomorphisms? Conversely,
can we see any symplectic homeomorphism (element of the $C^0$
closure of symplectic diffeomorphisms in the homeomorphisms in
general dimension, area-preserving homeomorphisms in dimension 2) as
an element of $\overline{\mathcal{H}}^{\gamma}$?

This last question is related to Oh's still open question whether
his group of "Hamiltonian homeomorphisms", called $Hameo$, equals or
not the group of symplectic homeomorphisms \cite{OH}.

\section{Application to the Hamilton-Jacobi equation.}\label{Hamilton Jacobi}

Let $H$ be a smooth Hamiltonian function on $\R\times\RN$. We
consider the evolution Hamilton-Jacobi equation $(HJ)$:
$$\derive{u}{t}+H\left(t,x,\derive{u}{x}\right)=0,$$
where $u:\R\times\R^n\to\R,\,(t,x)\mapsto u(t,x)$ satisfies an
initial condition $u(0,x)=u_0(x)$. First, we remind the reader of
the construction of a variational solution of $(HJ)$ (see for
example \cite{VX} or \cite{OV}).

\subsection{Recall on variational solutions of $(HJ)$.}

Let us denote by $\Lambda_0$ the graph of $du_0$ and call it the
initial submanifold. In fact, the following construction can be made
for any Lagrangian submanifold $\Lambda_0\subset\RN$. We consider
$\Sigma=\hat{H}^{-1}(\{0\})\subset\R^{2+2n}$. A geometric solution
of $(HJ)$ is a Lagrangian submanifold $L$ that satisfies
$\Lambda_0\leq L\leq\Sigma$. For example, the graph of the
differential of a smooth function $u$ is a geometric solution if and
only if $u$ itself is solution of $(HJ)$.

With the help of the flow $\phi_{\hat{H}}^t$, we can construct a geometrical solution $L_H=\bigcup_{t\in I}{\phi_{\hat{H}}^t(\Lambda_0)}$, where $I$ is an open interval containing $[0,1]$ and such that $\rho_{\alpha}=1$ on $I$. The Lagrangian submanifold $L_H$ obtained is an element of $\mathcal{L}(\R^{2+2n})$.

For any element $L\in\mathcal{L}(\R^{2k})$, we can associate a function $u_L$ on $\R^{2k}$ by the following method.

 Let $S:\R^k\times\R^q\to\R$ be a g.f.q.i of $L$. Denote by $1_z$ the fundamental class in $H^0({z})$. Then, we define $u_L$ by
$$u_L(z)=c(1_z,S|_{{z}\times\R^q}),$$
with notations of section \ref{invariants}. The function $u_L$ is everywhere $C^0$, and it is proved in \cite{OV}, that $u_L$ is $C^k$ on a dense open set, for $k\geq 1$. Moreover, when it is defined, we have $(x,du_L(x))\in L$.
Therefore, the function $u_{L_H}$ is a solution of $(HJ)$ on any open set on which it is smooth.

\medskip
We are now going to prove an interesting property of the elements of
$\overline{Ham}^{\gamma_2}$, which is the fact that we can extend to
them the construction of a variational solution of $(HJ)$.

\subsection{Extension to the completion}

\begin{prop}\label{propHJ} Let $H$ and $K$ be two Hamiltonian functions, and $u_{L_H}$, $u_{L_K}$ the solution obtained by the above method with the same initial submanifold $\Lambda_0$. Then,
$$\|u_{L_H}-u_{L_K}\|_{C^0}\leq {{\gamma}}_2(H,K).$$
\end{prop}

That leads us to the following definition.

\begin{definit} Let $H\in \overline{Ham}^{\gamma_2}$. A continuous function $u$ will be called \textnormal{generalized variational solution} of $(HJ)$ for $H$, if there exists a Cauchy sequence $(H_k)$ in $Ham$ representing $H$ and such that the sequence of solutions $(u_{L_{H_k}})$ $C^0$-converges to $u$.
\end{definit}

Therefore, proposition \ref{propHJ} implies the following statement:

\begin{theo} For each initial condition $u_0$, any element $H$ in the completion $\overline{Ham}^{\gamma_2}$ admits a unique generalized variational
solution $u_H$. Moreover, the so constructed map
$\overline{Ham}^{\gamma_2}\to C^0$ is continuous.
\end{theo}
 In particular, any Hamiltonian function in
$\mathfrak{F}^{\infty}$ (see definition \ref{la famille autonome})
admits a unique generalized variational solution.

\medskip
\noindent \textit{Proof. } Let $(H_k)\in Ham$ be a Cauchy sequence
for ${\gamma}_2$ representing an element $H\in
\overline{Ham}^{\gamma_2}$. Then, proposition \ref{propHJ} implies
that $(u_{L_{H_k}})$ is a Cauchy sequence in $C^0$ and hence
converges to a continuous function $u$. Moreover, if $(H_k)$ and
$(F_k)$ are two equivalent Cauchy sequences for ${\gamma}_2$, then
proposition \ref{propHJ} also implies that $(u_{L_{H_k}})$ and
$(u_{L_{F_k}})$ are equivalent, and hence converge to the same
limit. It gives the existence and the unicity.

The continuity of the map $\overline{Ham}^{\gamma_2}\to C^0$ is also
an immediate consequence of Proposition \ref{propHJ}. $\quad\Box$

\medskip
 To prove proposition \ref{propHJ},
we first prove the following lemma:

\begin{lem}\label{gamma et Hofer} For any $L\in\mathcal{L}$, we have
$$\|u_L\|_{C^0}\leq\gamma(L).$$
\end{lem}

\textit{Proof. }Since $L$ coincides with the zero section out of a compact set, $u_L$ has a compact support. It follows that $\|u_L\|_{C^0}\leq\max(u_L)-\min(u_L)$. We will prove that $\min(u_L)\geq c(1,L)$. It will also imply that $\max(u_L)\leq c(\mu,L)$ by Poincaré duality. Indeed, using the fact that $c(\mu,L)=-c(1,\overline{L})$ and that for all $z$, $\mu_z=1_z$, we have $u_L=-u_{\overline{L}}$.

Let $z\in\R^k$, and $S:\R^k\times\R^q\to\R$ be a g.f.q.i of $L\subset\R^{2k}$. Then, $S|_{\{z\}\times\R^q}$ is a g.f.q.i. of the reduction of $L$ by the coisotropic submanifold $\{z\}\times\R^k\subset\R^{2k}$. Therefore, by lemma \ref{red2}, we get $c(1_z,S|_{\{z\}\times\R^q})\geq c(1,S)$, for all $z$ and hence  $\min(u_L)\geq c(1,L)$ as required.$\quad\Box$

\medskip
\textit{Proof of proposition \ref{propHJ}. }The proposition comes from a sequence of inequalities:
$$\|u_{L_H}-u_{L_K}\|_{C^0}\leq\gamma(L_H,L_K)\leq\tilde{\gamma}(\phi_{\check{\hat{H}}},\phi_{\check{\hat{K}}})\leq {{\gamma}}_2(H,K).$$
The third inequality comes from the first inequality in proposition \ref{inegalites distances}. The second one is proved in \cite{CV}. Finally, the first one comes from the lemma \ref{gamma et Hofer} above and proposition 3.3 in \cite{V1}, which states that for all $u,v\in H^*(\R^n)$, $c(u\cup v,L_1+L_2)\leq c(u,L_1)+c(v,L_2)$, where $L_1+L_2=\{(q,p_1+p_2)\,|\,(q,p_1)\in L_1, (q,p_2)\in L_2 \}$.

 Indeed, for $u=v=1_{(t,x)}$, $L_1=(L_H-L_K)|_{(t,x)}$ and $L_2=L_K|_{(t,x)}$, we get $c(1_{(t,x)},L_H|_{(t,x)})-c(1_{(t,x)},L_K|_{(t,x)})\leq -c(1_{(t,x)},(L_H-L_K)|_{(t,x)})$. Then, lemma \ref{gamma et Hofer} gives  $-c(1_{(t,x)},(L_H-L_K)|_{(t,x)})\leq\gamma(L_H-L_K)=\gamma(L_H,L_K)$. By exchanging $H$ and $K$ and taking the supremum over $(t,x)$, we obtain $\|u_{L_H}-u_{L_K}\|_{C^0}\leq\gamma(L_H,L_K)$ as required.$\quad\Box$

\medskip
\noindent \textit{Remark and Question.} Joukovskaia proved in
\cite{J} that for Hamiltonian functions that are convex in $p$,
variational solutions of $(HJ)$ coincide with viscosity solutions
(These are a notion of weak solution introduced by Crandall and
Lions in \cite{CL} that has shown its efficiency in a lot of domains
of applications including optimal control and differential games,
front propagation problems, finance, image theory.... ). We are
tempted to use it together with some convergence result on viscosity
solutions, to prove that our generalized variational solution is a
viscosity solution. This would give another interpretation of our
notion of solution, and since our solution is continuous, it would
also give a continuity result on viscosity solutions.

However, since we developed our theory in the context of compactly
supported Hamiltonians, we cannot reason on Hamiltonian functions
convex in $p$. That leads us to our question : Can one define a
completion with similar properties for a class of Hamiltonian
functions convex in $p$?

\appendix
\section{Appendix: Proof of inequalities}

In this appendix we prove proposition \ref{inegalites distances} and lemma \ref{lemme capacité}. All those inequalities are based on the reduction inequality stated in proposition \ref{red1}.

\subsection{Inequality between $\tilde{\gamma}$ and $\gamma$.}\label{tilde}

We first prove the inequality $\gamma\geq\tilde{\gamma}$.

Let $\varphi$ be a Hamiltonian diffeomorphism, and $L\in\mathcal{L}$. We wish to show that $\gamma(\varphi(L)-L)\leq\gamma(\varphi)$. If we denote by $N$ the zero section of $\R^{2n}=T^*\R^n$, there exists a Hamiltonian isotopy $\psi^t$ such that $L=\psi^1(N)$. Therefore, we just need to prove $\gamma(\varphi(N))\leq\gamma(\varphi)$. Indeed, if we assume this inequality, then $\gamma(\varphi(L)-L)=\gamma(\varphi\circ\psi^1(N)-\psi^1(N))=\gamma(\psi^{-1}\circ\varphi\circ\psi^1(N)-N)$, using formula $(2.1)$ in \cite{CV}. Then, by assumption we get $\gamma(\varphi(L)-L)\leq \gamma(\psi^{-1}\circ\varphi\circ\psi^1)=\gamma(\varphi)$.

Let us prove now that $\gamma(\varphi(N))\leq\gamma(\varphi)$. We denote by $\Delta_p$ the diagonal in $\R^p\times\R^p$, and by $\Phi$ the symplectic identification $\overline{\R^{2n}}\times\R^{2n}\to T^*\Delta_{2n}$. Recall that $\widetilde{\Gamma_{\varphi}}$ is by definition the image of the graph $\Gamma_{\varphi}$ of $\varphi$. Clearly, $\varphi(N)$ is identified to the symplectic reduction of $N\times\Gamma_{\varphi}\subset\R^{6n}$ by the coisotropic linear subspace $\Delta_{2n}\times\R^{2n}$. It is therefore identified to the reduction of $N\times\widetilde{\Gamma_{\varphi}}$ by $W=(Id_{\RN}\times\Phi)(\Delta_{2n}\times\R^{2n})$. One can easily show that for all $L\in\mathcal{L}$, $\gamma(N\times L)=\gamma(L)$. In particular, $\gamma(\varphi)=\gamma(N\times\widetilde{\Gamma_{\varphi}})$, and the proof will be achieved if we prove the following proposition.

\begin{prop}[Reduction Inequality]\label{red1} For every Lagrangian submanifold $L$ in $\R^{2n}$ and every linear coisotropic subspace $W$ of $\R^{2n}$, we have $\gamma(L)\geq\gamma(L_W)$, where $L_W$ denotes the image of $L$ by reduction by $W$.
\end{prop}

We first prove the following lemma.

\begin{lem}\label{red2} Let $L$ be a Lagrangian submanifold in a cotangent bundle of the form $T^*M=T^*B\times\R^{2k}$. Consider the two coisotropic submanifolds $X=T^{\ast}B\times\{x_0\}\times\R^n$ and $Y=T^{\ast}B\times\R^n\times\{0\}$. Denote by $L_X$ and $L_Y$ the reductions of $L$ by respectively $X$ and $Y$. Then
 $$c(1,L_X)\geq c(1,L)=c(1,L_Y),$$
$$c(\mu_B,L_X)\leq c(\mu_M,L)=c(\mu_B,L_Y).$$
\end{lem}
\textit{Proof. } We start the proof by showing that $c(1,L_X)\geq
c(1,L)$. Let us fix $\lambda\in\R$ and  consider the inclusion
$i:B\simeq\{0\}\times B\to M$. Let $S$ be a g.f.q.i. of $L$ defined on a bundle $\pi:E\to M$. Then the function $S_X=S|_{\pi^{-1}(B\times\{x_0\})}$ is a generating function for $L_X$. Since $S_X$ is a restriction of $S$, we have an inclusion of the sublevels $S_X^{\lambda}\subset S^{\lambda}$, which induces a morphism $i_{\lambda}:H^{\ast}(S^{\lambda},S^{-\infty})\to H^{\ast}(S_X^{\lambda},S_X^{-\infty})$. The naturality of Thom isomorphism and the fact that all different inclusions commute make the following diagram commutative.
\begin{equation*}
\begin{CD}
H^{\ast}(B) @>T>> H^{\ast}(S_X^{\infty},S_X^{-\infty}) @>{j_{X,\lambda}^{\ast}}>> H^{\ast}(S_X^{\lambda},S_X^{-\infty})\\
@A{i^{\ast}}AA  @AA{i_{\infty}}A  @AA{i_{\lambda}}A\\
H^{\ast}(M) @>T>> H^{\ast}(S^{\infty},S^{-\infty}) @>{j_{\lambda}^{\ast}}>> H^{\ast}(S^{\lambda},S^{-\infty})
\end{CD}
\end{equation*}
Suppose now that $j_{X,\lambda}^{\ast}\circ T(1)\neq 0$. Then
$i_{\lambda}\circ j_{\lambda}^{\ast}\circ T(1)=j_{X,\lambda}^{\ast}\circ T\circ i^{\ast}(1)=j_{X,\lambda}^{\ast}\circ T(1)\neq 0$ hence $j_{\lambda}^{\ast}\circ T(1)\neq 0$. That proves $c(1,L_X)\geq c(1,L)$.

In the case of $L_Y$, we also have an explicit generating function, constructed as follows. Since $\R^k$ is contractible we can suppose that the fibers of $\pi$ do not depend on the second coordinate of $M$. Denote by $i:B\simeq
B\times\{0\}\to E$ the inclusion and by $\tau:B\times\R ^k\to B$ the trivial bundle of rank $k$ over $B$. Consider the vector bundle over $B$, $\rho=\tau\oplus i^{\ast}\pi$ whose total space is $F=\pi^{-1}( B\times\{0\}) \times \R^n$. Then, the function $S_Y$, defined for all $v\in B$ and $(x,\xi)\in \rho^{-1}(v)$ by $S_Y(v;x,\xi)=S(v,x;\xi)$ is a g.f.q.i for $L_Y$.
The map $f:E\to F,\ (v,x;\xi)\mapsto(v;x,\xi)$ is a diffeomorphism and satisfies $S_Y\circ f=S$. Therefore, we have $S_Y^{\lambda}=f(S^{\lambda})$, an isomorphism $H^{\ast}(S^{\lambda},S^{-\infty})\simeq H^{\ast}(S_Y^{\lambda},S_Y^{-\infty})$ and a commutative diagram
\begin{equation*}
\begin{CD}
H^{\ast}(B) @>T>> H^{\ast}(S_Y^{\infty},S_Y^{-\infty}) @>{j_{Y,\lambda}^{\ast}}>> H^{\ast}(S_Y^{\lambda},S_Y^{-\infty})\\
@A{i^{\ast}}AA  @AA{\simeq}A  @AA{\simeq}A\\
H^{\ast}(M) @>T>> H^{\ast}(S^{\infty},S^{-\infty}) @>{j_{\lambda}^{\ast}}>> H^{\ast}(S^{\lambda},S^{-\infty})
\end{CD}
\end{equation*}
The previous argument gives $c(1,L_Y)\geq c(1,L)$. The reverse inequality is obtained from the same diagram with the inclusion $i$ replaced by the projection $p:M\to B$ (which reverses vertical arrows).

Finally, $c(\mu_B,L_N)\leq c(\mu_M,L)=c(\mu_B,L_Y)$ is obtained from  $c(1,L_X)\geq c(1,L)=c(1,L_Y)$ by Poincar{\'e} duality, by noticing that $\overline{L_X}=\overline{L}_X$ and $\overline{L_Y}=\overline{L}_Y$.$\quad\Box$

\begin{lem}\label{decomp} Let $W$ be a coisotropic linear subspace of $\R^{2n}$. Denote by $N$ the zero section of $\R^{2n}=T^*\R^n$. Then there exists a decomposition in linear isotropic subspaces $\R^{2n}=N_1\oplus V_1\oplus N_2\oplus V_2\oplus N_3\oplus V_3$, where $N=N_1\oplus N_2 \oplus N_3$ and each $N_i\oplus V_i$, $i=1,2,3$ is a symplectic subspace, such that $W=N_1\oplus V_1\oplus N_2\oplus V_3$.
\end{lem}

\textit{Proof. }
Let us first recall that if $W$ is coisotropic with symplectic orthogonal $W^{\omega}\subset W$, any subspace $F$ such that $F\oplus W^{\omega}=W$ is symplectic. Indeed, since $F\subset W$, $F\cap F^{\omega} = F\cap F^{\omega}\cap W = F \cap(F\oplus W^{\omega})^{\omega}= F\cap W^{\omega}= \{0\}$.

If there exists a decomposition as in the lemma, then $W^{\omega}=N_2\oplus V_3$. Therefore we set $N_2=W^{\omega}\cap N$. Then, we define $N_1$ as one complementary of $N_2$ in $W\cap N$, and $F_1$ as one complementary of $W^{\omega}$ in $W$, containing $N_1$. By the above remark, $F_1$ is symplectic, and we can choose $V_1$ as one Lagrangian complementary of $N_1$ in $F_1$.

Then, we define $V_3$ as a complementary of $N_2$ in $W^{\omega}$. Since $W\cap N=N_1\oplus N_2$, $V_3\cap N={0}$, and we can define $N_3$ as a complementary of $N_1\oplus N_2$ in $N$. Then, $F_3=N_3\oplus V_3$ is symplectic since it is a complementary of $(N_1\oplus N_2\oplus F_3)^{\omega}$ in $N_1\oplus N_2\oplus F_3$.

Finally, we define $F_2$ as a complementary of $F_1\oplus F_3$ in $\R^{2n}$. Then, $F_2$ is symplectic for a similar reason as $F_3$, and we can define $V_2$ as a Lagrangian complementary of $N_2$ in $F_2$. The decomposition $\R^{2n}=N_1\oplus V_1\oplus N_2\oplus V_2\oplus N_3\oplus V_3$ satisfies all the requirements of lemma \ref{decomp}.$\quad\Box$

\medskip
\textit{Proof of proposition \ref{red1}.}
Since the linear symplectic group acts transitively on the set of all pairs of complementary Lagrangian subspaces (see proposition 7.4 in Chapter 1 of
\cite{LM}), and since the space of Lagrangian subspaces which
are complementary to the zero section $N$ is path connected, there
exists a symplectic isotopy $\Psi^t$ of $\RN$ such that $\Psi^0=Id$ and that
$\Psi^1$ lets all the elements of $N$ invariant and maps $V$ on
$V_1\oplus V_2\oplus V_3$.
Since $\RN$ is simply connected, that isotopy is Hamiltonian.

 The reduction of $L$ by $W$ is identified with the reduction of
$\Psi^1(L)$ by $\Psi^1(W)$. Therefore, applying twice the lemma
\ref{decomp}, we get $\gamma(L_W)\leq\gamma(\Psi^1(L))$.
But, by proposition 2.6 in \cite{V1}, we have
$\gamma(L)=\gamma(\Psi^1(L))$. That concludes the proof of proposition \ref{red1}.$\quad\Box$

\medskip
\noindent
\textit{Remark. } Note that in the end of the previous proof, lemma \ref{decomp} also implies $c(1,L_W)\geq c(1,L)$. That will be useful in the proof of lemma \ref{lemme capacité}.

\subsection{Inequalities involving the ``suspended distances''.}\label{suspended}

We now prove the inequality $\gamma_u(H,K)\leq\hat{\gamma}(H,K)$, for any $H$,$K$ Hamiltonian functions. It is sufficient to prove that for all Hamiltonian functions $H$,$K$, all $s$ in $[0,1]$, and all $\alpha$ large enough, $\gamma(\phi_K^{-s}\phi_H^{s})\leq\gamma(\phi_{\hat{K}}^{-s}\phi_{\hat{H}}^{s})$. We will prove that the graph of $\phi_K^{-s}\phi_H^{s}$ can be obtained by reduction of the graph of $\phi_{\hat{K}}^{-s}\phi_{\hat{H}}^{s}$, and then use proposition \ref{red1}.

We denote by $\hat{\Phi}^s$ the flow at time $s$ of the Hamiltonian $\hat{H}:(s;t,\tau,x)\mapsto \rho_{\alpha}(\tau)\tau+\rho_{\alpha}(t)H(t;x)$. By direct computation, we get
$$\hat{\Phi}^s(t,\tau,x)=(t(s),\tau(s),x(s)),$$ with
\begin{eqnarray*}
t(s) & = & t+\int_0^s (\rho_{\alpha} '(\tau(\sigma))\tau(\sigma)+\rho_{\alpha}(\tau(\sigma))d\sigma\\
\tau(s) & = & \tau-\int_0^s (\rho_{\alpha} '(t(\sigma))H(t(\sigma),x(\sigma)+\rho_{\alpha}(t(\sigma))\derive{H}{t}(t(\sigma),x(\sigma)))d\sigma\\
\end{eqnarray*}
and $x(s)$ solution of $\dot{x}(s)=\rho_{\alpha}(t(s))X_H(t(s),x(s))$.
If we denote $M=\max(\|\rho_{\alpha}\|_{C^1},\|H\|_{C^1})$, we see that $\tau(s)\in [\tau-|s|M^2,\tau+|s|M^2]$. Suppose $\tau\in [-M^2-2M,M^2+2M]$ and ${\alpha}$ is large enough, then $\rho_{\alpha}(\tau(s))=1$ and $t(s)=t+s$. Hence $x(s)=(\phi_H)_t^{t+s}(x)$.
We set
\begin{eqnarray*}
I_H(s,t,x) & = & -\int_0^s (\rho_{\alpha} '(t(\sigma))H(t(\sigma),x(\sigma)+\rho_{\alpha}(t(\sigma))\derive{H}{t}(t(\sigma),x(\sigma)))d\sigma\\
 & = & H(t,x)-H(t+s,\phi_t^{t+s}(x)),
\end{eqnarray*}
and $J(s,t,x)=I_H(s,t,x)+I_K(-s,t+s,(\phi_H)_t^{t+s}(x))$.
Then, we can write the expression of the composition:
$$\phi_{\hat{K}}^{-s}\phi_{\hat{H}}^{s}(t,\tau,x)=(t,\tau+J(s,t,x),(\phi_K)_t^{t-s}(\phi_H)_t^{t+s}(x)).$$

We can now compute the intersection of the graph $\Gamma_{\phi_{\hat{K}}^{-s}\phi_{\hat{H}}^{s}}$ with the set $U=[-1,1]\times\R\times[-M^2,M^2]\times\R\times\RN\times\RN$, and its image by the natural identification $\Psi:\R^{4+4n}\to T^*\Delta_{2+2n}$. We get
$$\widetilde{\Gamma}_{\phi_{\hat{K}}^{-s}\phi_{\hat{H}}^{s}}\cap\Psi(U)=\{(t,J(s,t,x),\tau+\frac{1}{2}J(s,t,x),0,z(x)) \,|$$
$$\qquad\qquad\qquad (t,\tau,x)\in[0,1]\times[-M^2,M^2]\times\RN,\ z(x)\in\widetilde{\Gamma}_{(\phi_K)_t^{t-s}(\phi_H)_t^{t+s}} \}.$$

Consider the coisotropic submanifold $W=\{0\}\times\R\times\{0\}\times\R\times\R^{4n}$. Since $\tau+\frac{1}{2}J(s,t,x)=0$ implies $\tau\in[-M^2-2M,M^2+2M]$, and since $W\subset \Psi(U)$, we see that $\widetilde{\Gamma}_{\phi_K^{-s}\phi_H^{s}}$ is obtained from $\widetilde{\Gamma}_{\phi_{\hat{K}}^{-s}\phi_{\hat{H}}^{s}}$ by reduction by $W$. By proposition \ref{red1}, we get $\gamma(\widetilde{\Gamma}_{\phi_K^{-s}\phi_H^{s}})\leq\gamma(\widetilde{\Gamma}_{\phi_{\hat{K}}^{-s}\phi_{\hat{H}}^{s}})$ and hence $\gamma(\phi_K^{-s}\phi_H^{s})\leq\gamma(\phi_{\hat{K}}^{-s}\phi_{\hat{H}}^{s})$.$\quad\Box$

\bigskip
We are now going to prove $\gamma_u(H,K)\leq\check{\gamma}(H,K)$. The idea of the proof is the same as the previous one: we show that for any $s\in[0,1]$, $\widetilde{\Gamma}_{\phi_K^{-s}\phi_H^{s}}$ is obtained by reduction of $\widetilde{\Gamma}_{{\phi_{\check{K}}}^{-1}\phi_{\check{H}}}$, for $\alpha$ large enough.

Recall that by definition, $\check{H}(s;t,\tau,x)=\rho_{\alpha}(t)tH(st;x)$. As above, we compute the flow : $\phi_{\check{H}}^s(t,\tau,x)=(t(s),\tau(s),x(s))$, and we obtain
\begin{eqnarray*}
t(s) & = & t \\
\tau(s) & = & \tau + I_H(s,t,x)
\end{eqnarray*}
where $I_H(s,t,x)=\rho_{\alpha}(t)sH(st,x(s))-\rho_{\alpha}'(t)t\int_0^s{H(\sigma t,x(\sigma))d\sigma}$ and $x(s)$ is solution of $\dot{x}(s)=\rho_{\alpha}(t)t X_H(st,x(s))$. For $t\in[-1,1]$ and $\alpha\geq 1$, it gives $x(s)=\phi^{ts}(x)$.

Similarly as above, we set $J(s,t,x)=I_H(s,t,x)+I_K(-s,t+s,(\phi_H)^{ts}(x))$, the set $U=[-1,1]\times\R\times\R^2\times\RN\times\RN$ and the identification $\Psi:\R^{4+4n}\to T^*\Delta_{2+2n}$. The graph can be written this way:
$$\widetilde{\Gamma}_{\phi_{\check{K}}^{-s}\phi_{\check{H}}^{s}}\cap\Psi(U)=\{(t,J(s,t,x),\tau+\frac{1}{2}J(s,t,x),0,z(x)) \,|$$
$$\qquad\qquad\qquad (t,\tau,x)\in[0,1]\times\R\times\RN,\ z(x)\in\widetilde{\Gamma}_{\phi_K^{-st}\phi_H^{st}} \}.$$
Now, we see that  $\widetilde{\Gamma}_{\phi_K^{-t}\phi_H^{t}}$ is the reduction of $\widetilde{\Gamma}_{{\phi_{\check{K}}}^{-1}\phi_{\check{H}}}$ by the coisotropic manifold $W=\{t\}\times\R\times\{0\}\times\R\times\R^{4n}$. Using lemma \ref{red2} twice, we conclude that for all $t\in[0,1]$, $\gamma(\phi_K^{-t}\phi_H^{t})\leq\gamma(\phi_{\check{K}}^{-1}\phi_{\check{H}})$.$\quad\Box$

\subsection{Proof of lemma \ref{lemme capacité}.}
It is sufficient to show that $c(V)\leq c(\R^2\times V)$ for all
open subset $V\in\RN$. Let $H$ be a Hamiltonian function with
support in $V$. We just have to find a Hamiltonian function $K$ with
support in $V\times\R^2$ satisfying the inequality $c_+(H)\leq
c_+(K)$. If we set $K=\check{H}_{\alpha}$ for $\alpha$ large enough,
$K$ has support in $\R^2\times V$, and we saw in particular in the
previous proof that $\widetilde{\Gamma}_{\phi_H^{1}}$ is the
reduction of $\widetilde{\Gamma}_{\phi_{\check{H}}}$. Therefore, by
the remark that ends section \ref{tilde}, we have $c_+(H)\leq
c_+(K)$ as required.$\quad\Box$

\nocite{*}
\bibliographystyle{plain}
\bibliography{article17}

\end{document}